\documentclass{elsarticle}
\usepackage{lineno,hyperref}

\usepackage{amsmath,amsfonts,amsthm}
\usepackage{fixmath}
\usepackage{amssymb,latexsym}
\usepackage{mathrsfs,amsbsy}
\usepackage{dsfont}
\usepackage{algorithm}
\usepackage{listings}
\usepackage{eucal}
\usepackage{colortbl}
\usepackage{multirow}
\usepackage{float}
\graphicspath{{./Figure_submit/}}
\usepackage{booktabs}
\usepackage{caption}
\usepackage{subfigure}
\usepackage{microtype}
\usepackage{enumerate}
\usepackage{hyperref}
\usepackage{afterpage}
\usepackage{xcolor}

\usepackage{setspace}
\usepackage[margin=1.5in]{geometry}
\usepackage{tikz} 
\usetikzlibrary{shapes,arrows} 

\modulolinenumbers[5]

\makeatletter
\def\ps@pprintTitle{%
	\let\@oddhead\@empty
	\let\@evenhead\@empty
	\def\@oddfoot{\footnotesize\itshape
		{Submitted preprint} \hfill\today}%
	\let\@evenfoot\@oddfoot
}
\makeatother

\begin{document}

\begin{frontmatter}

\title{Construction of Diffeomorphisms with Prescribed Jacobian Determinant and Curl}



\author[mymainaddress]{Zicong Zhou\corref{mycorrespondingauthor}}
\cortext[mycorrespondingauthor]{Corresponding author}
\ead{zicongzhou818@sjtu.edu.cn}

\author[mysecondaryaddress]{Guojun Liao}
\ead{liao@uta.edu}

\address[mymainaddress]{Institute of Natural Sciences, Shanghai Jiao Tong University, 800 Dongchuan Rd, Shanghai 200240, China}

\address[mysecondaryaddress]{Department of Mathematics, University of Texas at Arlington, 701 S. Nedderman Dr, Arlington, TX 76019, USA}

\begin{abstract}
The variational principle (VP) is designed to generate non-folding grids (diffeomorphisms) with prescribed Jacobian determinant (JD) and curl. Its solution pool of the original VP is based on an additive formulation and, consequently, is not invariant in the diffeomorphic Lie algebra. The original VP works well when the prescribed pair of JD and curl is calculated from a diffeomorphism, but not necessarily when the prescribed JD and curl are not known to come from a diffeomorphism. This issue is referred as the mismatched pair problem. In spite of that, the original VP works effectively in 2D grid generations. To resolve these issues, in this paper, we describe a new version of VP (revised VP), which is based on composition of transformations and, therefore, is invariant in the Lie algebra. The revised VP seems have overcome the inaccuracy of original VP in 3D grid generations. In the following sections, the mathematical derivations are presented. It is shown that the revised VP can calculate the inverse transformation of a known diffeomorphism. Its inverse consistency and transitivity of transformations are also demonstrated numerically. Moreover, a computational strategy is formulated based on the new version of VP to handle the mismatch issue and is demonstrated with preliminary result. 
\end{abstract}

\begin{keyword}
\texttt{adaptive grid generation}\sep \texttt{computational diffeomorphism}\sep \texttt{Jacobian determinant}\sep \texttt{curl}
\MSC[2010] 49Q10 \sep 49Q20 \sep 65K10 \sep 68W25 \sep 93B27 \sep 93B40
\end{keyword}

\end{frontmatter}


\section{Introduction}
Computational construction of diffeomorphism is an active research field in computational geometry. For instance, conformal differential geometry \cite{Gu} led by Gu achieved remarkable success in surface differential geometry. Whereas, the aims of this study is about how to characterize and control a meaningful distribution of grid points over a volumetric domain, such as in \cite{Grajewski,Liseikin}. This is the problem of adaptive generation of non-folding grids. One approach to the task is to find a differentiable and invertible transformation, i.e., a diffeomorphism, $\pmb{T}$ by controlling its JD which models local cell-size, such as in \cite{Brackbill, DacMos, Moser, Huang}. Its idea has been widely transfered and applied in constructing deformable image registration methods, such as in \cite{ChenY, Lee, Joshi, Joshi2, Sotiras}. A grid generation method in \cite{Cai}, the deformation constructs $\pmb{T}$ with prescribed JD, by a scalar $0<f \in C^1$, whose key component is the solution to a $\bf{divergence-curl}$ system (similar to the constraints of (\ref{ssdregular})). Its divergence is approximated by $f-1$ and curl is assigned by $\pmb{0}$, due to the challenge of realizing curl before hand. Consequently, grids generated by the deformation method is not unique due to the lack of curl information, in \cite{Liao}, which models the local cell-rotation. To overcome this problem, the original VP was proposed in \cite{ChenXi} and studied further in \cite{Zhou}. In 3D grid generations, the original VP only provides inaccurate approximations around the true solution. The authors attempted to analyze the uniqueness of such transformations, but it turned out cannot be completed for the reason that JD of a transformation is merely an approximation to the divergence of the transformation. Fortunately, this does not undermine VP to produce well approximated grids and in turn allows mismatched values of prescribed JD and curl under certain range. 

However, a novel image atlas construction method proposed in \cite{Zhou} utilizes the original VP and requires it (i) to accept wider ranges of mismatched JD and curl and (ii) to optimize JD and curl in a separable manner so one may investigate how each of them affect a diffeomorphism. Another limitation of the original VP is its consideration of small deformations $\pmb{T}=\pmb{id}+\pmb{u}\in H^{2}_{0}(\mathrm{\Omega})$ where $\pmb{u}$ is the displacement field and $\pmb{id}$ is the identity map (uniform grid), such as in \cite{Joshi}, whose function composition ``$\circ$" is approximated by $\pmb{T}_2 \circ \pmb{T}_1\approx\pmb{T}_{1}+\pmb{u}_{2}=\pmb{id}+\pmb{u}_{1}+\pmb{u}_{2}=\pmb{id}+\pmb{u}_{2}+\pmb{u}_{1}=\pmb{T}_{2}+\pmb{u}_{1}\approx\pmb{T}_1 \circ \pmb{T}_2$. From a computational perspective, it risks having function composition maps outside of the collection $\pmb{T}\in H^{2}_{0}(\mathrm{\Omega})$ that forms a diffeomorphism group \cite{Bauer, Joshi2}. Therefore, in this paper, to resolve (i), (ii) and to broaden the transformations that VP can characterize, we fundamentally revise VP to consider transformations taking composition as the left-translation. Surprisingly, this revision also overcome the inaccuracy of the original VP in 3D grid generation.

The structure of the paper is organized as follow. In section \ref{vpnew}, we reformulated VP to cope with composition as left-translation. In section \ref{demo}, examples are provided to demonstrate the effectiveness of the revised VP; In section \ref{demo2}, a computational strategy is proposed to handle the mismatch issue. The numerical experiments run with MatLab codes on a desktop PC with ADM Ryzen-9 12-core Processor, 16 GB RAM and NVIDA GeForce RTX 3080 GPU.

\section{New Version of Variational Principle}\label{vpnew}
Let a simply-connected, bounded $\mathrm{\Omega} \subset \mathbb{R}^{3}$ (similar in $\mathbb{R}^{2}$) be the domain and  $\pmb{x}=(x,y,z)\in\mathrm{\Omega}$. Let a scalar function $f_o>0$ and a vector-valued function $\pmb{g}_o$ on $\mathrm{\Omega}$ satisfy $\int_{\mathrm{\Omega}} f_o(\pmb{x})d\pmb{x} = |\mathrm{\Omega}|$ and $\nabla \cdot  \pmb{g}_o= 0$, respectively. Given $\pmb{\phi}_{o}\in H^{2}_{0}(\mathrm{\Omega})$, we look for a diffeomorphic transformation $\pmb{\phi}=\pmb{\phi}_{\pmb{m}} \circ \pmb{\phi}_{o}=\pmb{\phi}_{\pmb{m}}(\pmb{\phi}_{o})\in H^{2}_{0}(\mathrm{\Omega})$, where $\pmb{\phi}_{\pmb{m}}=\pmb{id}+\pmb{u}$ is an intermediate transformation that left-translates $\pmb{\phi}_{o}$ to $\pmb{\phi}$ and  implies $\delta\pmb{\phi}_{\pmb{m}}=\delta\pmb{u}$, that the cost functional --- sum of squared differences ($SSD$) is minimized:
\begin{equation}\label{ssdvar}
	SSD(\pmb{\phi}) = \frac{1}{2}\int_{\mathrm{\Omega}} [(\text{det}\nabla\pmb{\phi} - f_o)^2+|\nabla \times\pmb{\phi}- \pmb{g}_o|^2] d\pmb{x}
\end{equation}	
	
\begin{equation}\label{ssdregular}	
	\text{ subjects to }
	\left\{
	\begin{aligned}
		&
		\begin{aligned}
			\nabla \cdot  \pmb{u}& = f-1 \\
			\nabla \times \pmb{u}& = \pmb{g}
		\end{aligned}
	\end{aligned}\right.
\hspace{-0.5cm}
	\Rightarrow
		\mathrm{\Delta} \pmb{u} =\nabla f-\nabla \times \pmb{g} = \pmb{F} \text{ in } \mathrm{\Omega},
\end{equation}	
with  $\pmb{u} = \pmb{0}$ on $\partial\mathrm{\Omega}$ where $f$, $\pmb{g}$ and $\pmb{F}$ are control functions. Its variational gradient with respect to the control function $\pmb{F}$  can be derive as follows. Denote $P=\text{det}\nabla\pmb{\phi} - f_o$ and $\pmb{Q}=\nabla \times\pmb{\phi}- \pmb{g}_o$, then, for all $\delta\pmb{F}$ vanishing on $\partial\mathrm{\Omega}$,
\begin{equation*}
	\begin{aligned}
		\delta SSD(\pmb{\phi}_{\pmb{m}} &\circ \pmb{\phi}_{o})	=\int_{\mathrm{\Omega}} [(\text{det}\nabla\pmb{\phi} - f_o)\delta\text{det}\nabla(\pmb{\phi}_{\pmb{m}} \circ \pmb{\phi}_{o}) +(\nabla \times\pmb{\phi}- \pmb{g}_o)\cdot\delta\nabla \times(\pmb{\phi}_{\pmb{m}} \circ \pmb{\phi}_{o})]d\pmb{x}\\
		=&\int_{\mathrm{\Omega}} [P\delta\text{det}\nabla\pmb{\phi}_{\pmb{m}}\text{det}\nabla\pmb{\phi}_{o} +\pmb{Q}\cdot\delta
		\begin{pmatrix}
			\nabla\phi_{\pmb{m}3}\cdot(\pmb{\phi}_{o})_{y}-\nabla\phi_{\pmb{m}2}\cdot (\pmb{\phi}_{o})_{z}\\ 
			-\nabla\phi_{\pmb{m}3}\cdot(\pmb{\phi}_{o})_{x}+\nabla\phi_{\pmb{m}1}\cdot(\pmb{\phi}_{o})_{z} \\ 
			\nabla\phi_{\pmb{m}2}\cdot(\pmb{\phi}_{o})_{x}-\nabla\phi_{\pmb{m}1}\cdot(\pmb{\phi}_{o})_{y}    
		\end{pmatrix}]d\pmb{x}\\
		=\int_{\mathrm{\Omega}} [P&\text{det}\nabla\pmb{\phi}_{o}(\delta u_{1x}\phi_{\pmb{m}2y}\phi_{\pmb{m}3z}+\phi_{\pmb{m}1x}\delta u_{2y}\phi_{\pmb{m}3z}+\phi_{\pmb{m}1x}\phi_{\pmb{m}2y}\delta u_{3z}
		-\delta u_{1x}\phi_{\pmb{m}2z}\phi_{\pmb{m}3y}\\
		-\phi_{\pmb{m}1x}&\delta u_{2z}\phi_{\pmb{m}3y}-\phi_{\pmb{m}1x}\phi_{\pmb{m}2z}\delta u_{3y}-\delta u_{1y}\phi_{\pmb{m}3z}\phi_{\pmb{m}2x}-\phi_{\pmb{m}1y}\delta u_{3z}\phi_{\pmb{m}2x}-\phi_{\pmb{m}1y}\phi_{\pmb{m}3z}\delta u_{2x}\\
		+\delta u_{1y}&\phi_{\pmb{m}3x}\phi_{\pmb{m}2z}+\phi_{\pmb{m}1y}\delta u_{3x}\phi_{\pmb{m}2z}+\phi_{\pmb{m}1y}\phi_{\pmb{m}3x}\delta u_{2z}
		+\delta u_{1z}\phi_{\pmb{m}2x}\phi_{\pmb{m}3y}+\phi_{\pmb{m}1z}\delta u_{2x}\phi_{\pmb{m}3y}\\
		+\phi_{\pmb{m}1z}&\phi_{\pmb{m}2x}\delta u_{3y}
		-\delta u_{1z}\phi_{\pmb{m}2y}\phi_{\pmb{m}3x}-\phi_{\pmb{m}1z}\delta u_{2y}\phi_{\pmb{m}3x}-\phi_{\pmb{m}1z}\phi_{\pmb{m}2y}\delta u_{3x})\\
		+\begin{pmatrix}
			Q_1 \\ 
			Q_2 \\ 
			Q_3 \notag   
		\end{pmatrix}&\cdot
		\begin{pmatrix}
			\delta u_{3x}\phi_{o1y}+\delta u_{3y}\phi_{o2y}+\delta u_{3z}\phi_{o3y}-\delta u_{2x}\phi_{o1z}-\delta u_{2y}\phi_{o2z}-\delta u_{2z}\phi_{o3z}\\ 
			-\delta u_{3x}\phi_{o1x}-\delta u_{3y}\phi_{o2x}-\delta u_{3z}\phi_{o3x}+\delta u_{1x}\phi_{o1z}+\delta u_{1y}\phi_{o2z}+\delta u_{1z}\phi_{o3z}\\ 
			\delta u_{2x}\phi_{o1x}+\delta u_{2y}\phi_{o2x}+\delta u_{2z}\phi_{o3x}-\delta u_{1x}\phi_{o1y}-\delta u_{1y}\phi_{o2y}-\delta u_{1z}\phi_{o3y}    
		\end{pmatrix}]d\pmb{x}\\
		=\int_{\mathrm{\Omega}} [(P&\text{det}\nabla\pmb{\phi}_{o}\begin{pmatrix}
			\phi_{\pmb{m}2y}\phi_{\pmb{m}3z}-\phi_{\pmb{m}3y}\phi_{\pmb{m}2z} \\ 
			\phi_{\pmb{m}3x}\phi_{\pmb{m}2z}-\phi_{\pmb{m}2x}\phi_{\pmb{m}3z} \\ 
			\phi_{\pmb{m}2x}\phi_{\pmb{m}3y}-\phi_{\pmb{m}2y}\phi_{\pmb{m}3x} \notag   
		\end{pmatrix}
		+\begin{pmatrix}
			Q_2\phi_{o1z}-Q_3\phi_{o1y} \\ 
			Q_2\phi_{o2z}-Q_3\phi_{o2y} \\ 
			Q_2\phi_{o3z}-Q_3\phi_{o3y} \notag   
		\end{pmatrix})\cdot  \nabla \delta u_1\\
		+&(P\text{det}\nabla\pmb{\phi}_{o}\begin{pmatrix}
			\phi_{\pmb{m}3y}\phi_{\pmb{m}1z}-\phi_{\pmb{m}1y}\phi_{\pmb{m}3z} \\ 
			\phi_{\pmb{m}1x}\phi_{\pmb{m}3z}-\phi_{\pmb{m}1z}\phi_{\pmb{m}3x} \\ 
			\phi_{\pmb{m}3x}\phi_{\pmb{m}1y}-\phi_{\pmb{m}1x}\phi_{\pmb{m}3y} \notag   
		\end{pmatrix}
%
		+\begin{pmatrix}
			-Q_1\phi_{o1z}+ Q_3\phi_{o1x} \\ 
			-Q_1\phi_{o2z}+ Q_3\phi_{o2x} \\ 
			-Q_1\phi_{o3z}+ Q_3\phi_{o3x}\notag   
		\end{pmatrix})\cdot  \nabla \delta u_2\\
		+&(P\text{det}\nabla\pmb{\phi}_{o}\begin{pmatrix}
			\phi_{\pmb{m}1y}\phi_{\pmb{m}2z}-\phi_{\pmb{m}2y}\phi_{\pmb{m}1z} \\ 
			\phi_{\pmb{m}2x}\phi_{\pmb{m}1z}-\phi_{\pmb{m}1x}\phi_{\pmb{m}2z} \\ 
			\phi_{\pmb{m}1x}\phi_{\pmb{m}2y}-\phi_{\pmb{m}2x}\phi_{\pmb{m}1y} \notag   
		\end{pmatrix}+\begin{pmatrix}
			Q_1\phi_{o1y}- Q_2\phi_{o1x} \\ 
			Q_1\phi_{o2y}- Q_2\phi_{o2x} \\ 
			Q_1\phi_{o3y}- Q_2\phi_{o3x} \notag   
		\end{pmatrix})\cdot \nabla \delta u_3]d\pmb{x}\\
		=&\int_{\mathrm{\Omega}} (\pmb{A}_1\cdot \nabla \delta u_1
		+\pmb{A}_2\cdot \nabla \delta u_2+\pmb{A}_3\cdot \nabla  \delta u_3)d\pmb{x}.
	\end{aligned}
\end{equation*}
Here, the ``big vector"s are now denoted as $\pmb{A}_{i}$, where $i=1,2,3$. By $Green$'s identities with fixed boundary condition and for some $B_{i}$ such that $\mathrm{\Delta} B_{i}=-\nabla \cdot \pmb{A}_{i}$ and $i=1,2,3$, then it can be carried to,
\begin{equation}\label{ssd2delF}
	\begin{aligned}
		\delta SSD=&\int_{\mathrm{\Omega}} (-\nabla\cdot \pmb{A}_1 \delta u_1
		-\nabla\cdot\pmb{A}_2 \delta u_2-\nabla\cdot\pmb{A}_3  \delta u_3)d\pmb{x}\\
		=&\int_{\mathrm{\Omega}} (\mathrm{\Delta} B_1  \delta u_1
		+\mathrm{\Delta} B_2 \delta u_2+\mathrm{\Delta} B_3  \delta u_3)d\pmb{x}\\
		=&\int_{\mathrm{\Omega}} ( B_1 \delta \mathrm{\Delta}u_1
		+B_2 \delta \mathrm{\Delta}u_2+ B_3  \delta \mathrm{\Delta}u_3)d\pmb{x}\\
		=&\int_{\mathrm{\Omega}} (\pmb{B} \cdot \delta \mathrm{\Delta} \pmb{u})d\pmb{x}\\
		=&\int_{\mathrm{\Omega}} (\pmb{B} \cdot \delta \pmb{F})d\pmb{x} \qquad
		\Rightarrow \qquad
		\frac{\partial SSD}{\partial  \pmb{F}}=\pmb{B}.
	\end{aligned}
\end{equation}
To give a more completed view, the variational gradients of $SSD$ with respect to $f$ and $\pmb{g}$ are included as well. For arbitrary $\delta f$ and $\delta\pmb{g}$, one may get $ \delta \mathrm{\Delta} \pmb{u} = \nabla \delta f-\nabla \times \delta \pmb{g}$, then, from (\ref{ssd2delF}), it can be derived, 
\begin{equation*}
	\begin{aligned}	
	\delta SSD = &\int_{\mathrm{\Omega}} \pmb{B} \cdot(\nabla \delta f-\nabla \times \delta \pmb{g})d\pmb{x}\\
	=&\int_{\mathrm{\Omega}} [\pmb{B} \cdot\nabla \delta f- B_{1}(\delta g_{3y}+ \delta g_{2z})+ B_{2}(\delta g_{3x}- \delta g_{1z})-B_{3} (\delta g_{2x}+ \delta g_{1y})]d\pmb{x}\\
		=&\int_{\mathrm{\Omega}} [\pmb{B} \cdot\nabla \delta f
		+ \begin{pmatrix}
			0 \\ 
			B_{3} \\ 
			-B_{2}    
		\end{pmatrix}\cdot\nabla \delta g_{1}
		+\begin{pmatrix}
			-B_{3} \\ 
			0 \\ 
			B_{1}   
		\end{pmatrix} \cdot\nabla\delta g_{2}
		+\begin{pmatrix}
			B_{2} \\ 
			-B_{1} \\ 
			0 
		\end{pmatrix}\cdot\nabla \delta g_{3}
		]d\pmb{x}\\
		=&\int_{\mathrm{\Omega}} [-\nabla\cdot \pmb{b} \delta f 
		- \begin{pmatrix}
			B_{3y} - B_{2z}  \\ 
			-B_{3x} + B_{1z} \\ 
			B_{2x} - B_{1y}  
		\end{pmatrix} \cdot\delta \pmb{g}
		]d\pmb{x}\\
		=&\int_{\mathrm{\Omega}}-\nabla\cdot\pmb{B}  \delta f d\pmb{x}+ \int_{\mathrm{\Omega}} -\nabla \times  \pmb{B}  \cdot \delta \pmb{g}d\pmb{x}
	\end{aligned}
	\end{equation*}
\begin{equation}\label{ssddelfg}		
		\hspace{-1.5cm}\Rightarrow \quad
		\frac{\partial SSD}{\partial f}=-\nabla\cdot\pmb{B} \quad  \text{ and } \quad \frac{\partial SSD}{\partial  \pmb{g}}=-\nabla \times  \pmb{B}. 
\end{equation}
A gradient descent based algorithm is provided below. $\triangle t$ is the step-size of gradients.  Major computational costs occur in solving $Poisson$ equations by a Fast Fourier Transform $Poisson$ solver, denoted as FFT. Define $ratio=SSD_{final}/SSD_{initial}$ and let $\pmb{F}_{k}$ be the control functions to be optimized.
\begin{algorithm}[H]\label{alg}
	{
		\caption{$[ \pmb{\phi}_{k}, \pmb{\phi}^{k}_{\pmb{m}}]  =  $ revisedVP($f_o$, $\pmb{g}_o$, $\pmb{\phi}_{o}$)}
		\hrule
		\vspace{-0.2cm}
		\begin{itemize}
			\item[$\bullet$] 1: set $\triangle t$, $t_{tol}$, $ratio$, $ratio_{tol}$, $k_{max}$, $t_{up}$, $t_{down}$, $better=1$, $k=0$; 
			\item[$\bullet$] 2: initialize $\pmb{F}_{k}=\pmb{0}$ ($=1$ when optimize along JD direction of (\ref{ssddelfg})), $\pmb{b}_{k}=\pmb{0}$, $\pmb{\phi}_{k}=\pmb{\phi}_{o}$;	
			\item[$\bullet$] 3: while $\triangle t>t_{tol}$ and $ratio>ratio_{tol}$ and $k<k_{max}$;		
			\begin{itemize}
				\item[$\bullet$] 4: if $better=1$
				\begin{itemize}
					\item[$\bullet$] 5: $k=k+1$
					\item[$\bullet$] 6: solve for $\pmb{B}_k$ from $\mathrm{\Delta} B_i=-\nabla \cdot \pmb{A}_i$ by FFT where $i=1,2,3$;		
				\end{itemize}
				\item[$\bullet$] 7: update $\pmb{F}_{k}=\pmb{F}_{k-1}-\triangle t * \pmb{b}_{k}$;
				\item[$\bullet$] 8: solve for $\pmb{u}_k$ from		$\mathrm{\Delta}\pmb{u}_k=\pmb{F}_{k}$ by FFT 
				\item[$\bullet$] 9: update $\pmb{\phi}_{k}=\pmb{\phi}^{k}_{\pmb{m}}( \pmb{\phi}_{o})$ by interpolation, where $\pmb{\phi}^{k}_{\pmb{m}}=\pmb{id}+\pmb{u}_k$;
				\item[$\bullet$] 10: compute $SSD$ and $ratio$;
				\item[$\bullet$] 11: if $SSD$ decrease,
				\begin{itemize}
					\item[$\bullet$] 12: $better=1$ 
					\item[$\bullet$] 13: $\triangle t=\triangle t * t_{up}$;
					\item[$\bullet$] 14: $\pmb{F}_{k-1}=\pmb{F}_{k}$;
				\end{itemize}
				\hspace{0.3cm} else
				\begin{itemize}
					\item[$\bullet$] 15: $better=0$ 
					\item[$\bullet$] 16: $\triangle t=\triangle t * t_{down}$.
					\vspace{-0.3cm}
				\end{itemize}
			\end{itemize}
		\end{itemize}	
	}
\end{algorithm}
\section{Numerical Examples of Revised VP}\label{demo}
E.g.\ref{eg34} is an example of a 3D grid reconstruction by the revised VP which demonstrates this work can be applied in 3D scenario. For more intuitive and clear visualizations, the rest of the examples are provided in 2D only.  E.g.\ref{eg1} shows the revised VP is capable of constructing inverse transformation; E.g.\ref{eg2} and \ref{eg3} confirm that the computational solutions of the revised VP satisfy inverse consistency and transitivity, which are expected in a diffeomorphism group. In all figures, ``a" vs ``b" means coloured in red  ``a" superposes on coloured in black ``b", except for displacement fields comparison of E.g.\ref{eg34}, which ``c" vs ``d" means coloured in blue  ``c" superposes on coloured in green ``d".

\subsection{$\mathbf{Example}$: 3D Grid Reconstruction Comparison between the Original and Revised VPs}\label{eg34}
\begin{figure}[H]
	\centering
	\subfigure[$\pmb{\Phi}$]{\includegraphics[width=3.4cm,height=3.4cm]{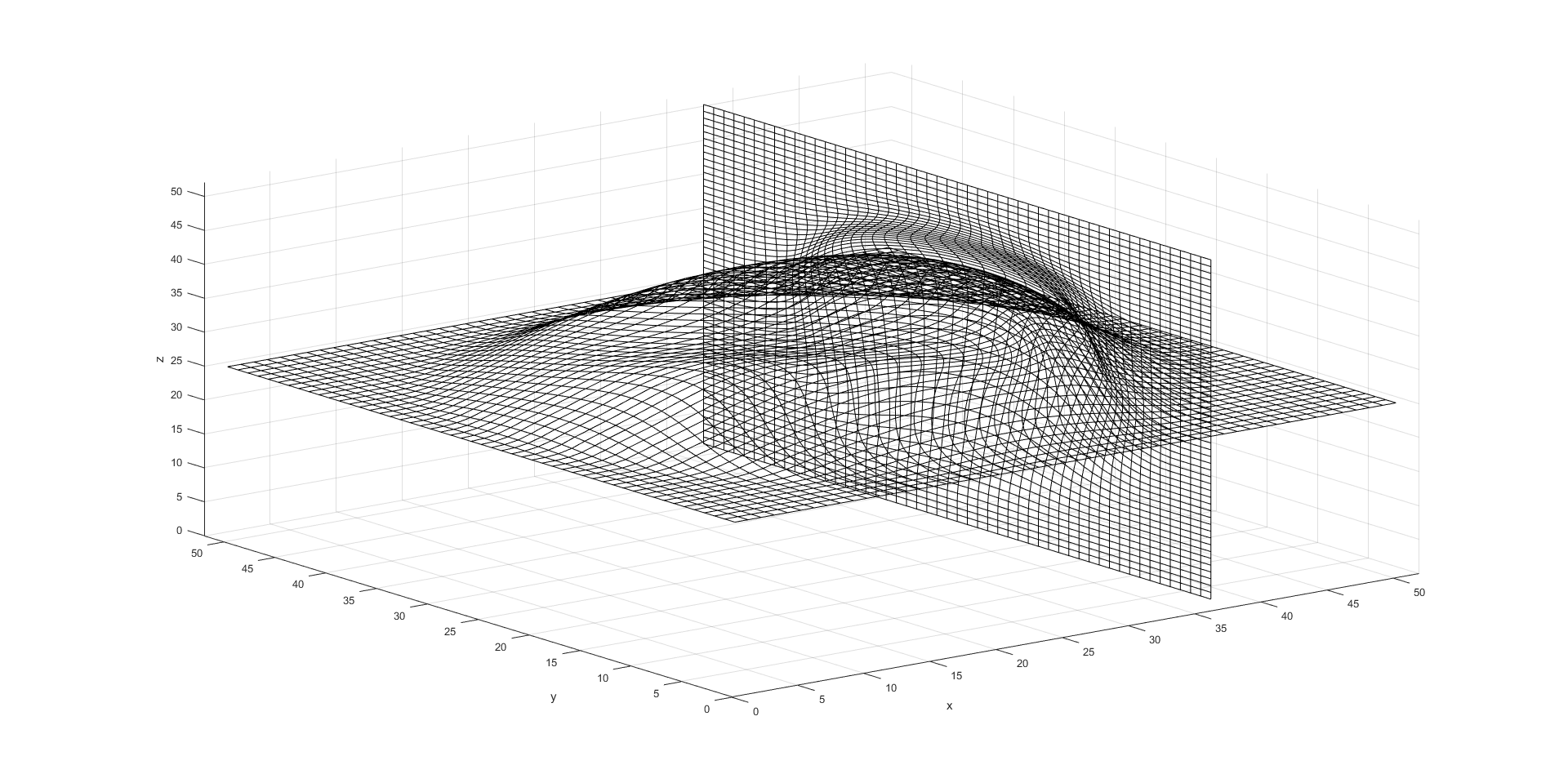}}
	\subfigure[$\pmb{\Phi}$ ($xz$-plain)]{\includegraphics[width=3.4cm,height=3.4cm]{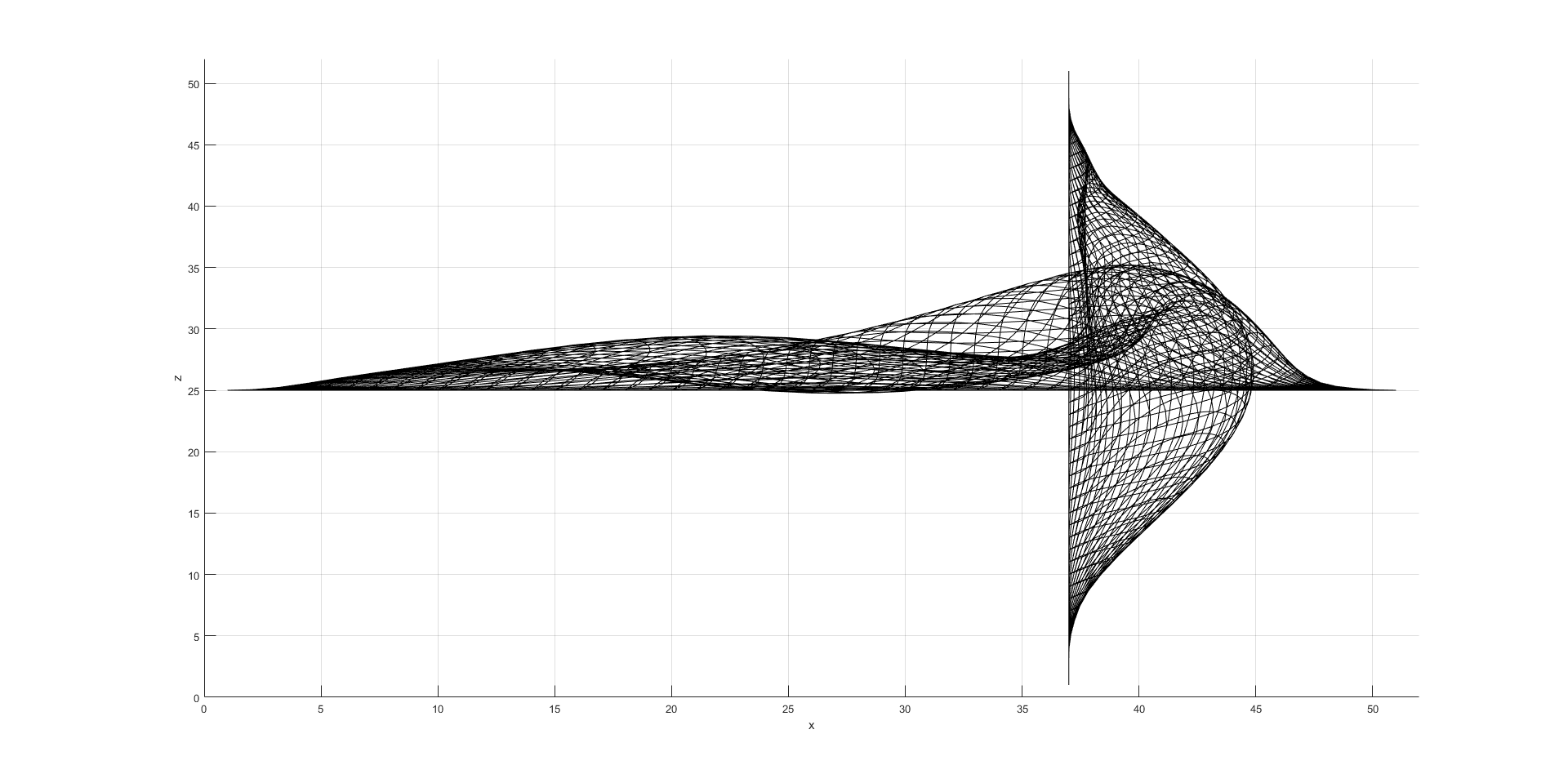}}
	\subfigure[$\pmb{U}$]{\includegraphics[width=3.4cm,height=3.4cm]{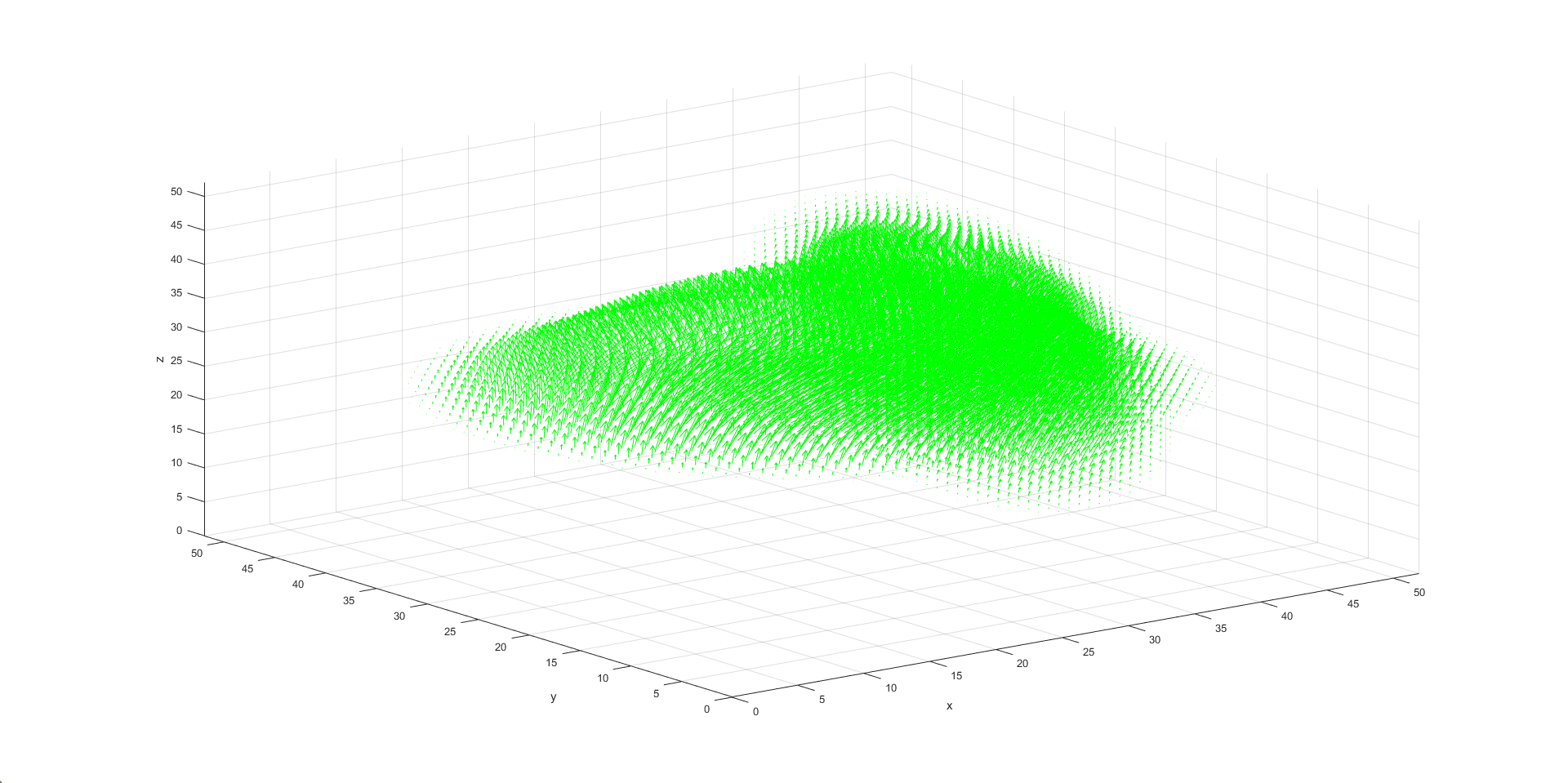}}
	\subfigure[$\pmb{U}$ ($xz$-plain)]{\includegraphics[width=3.4cm,height=3.4cm]{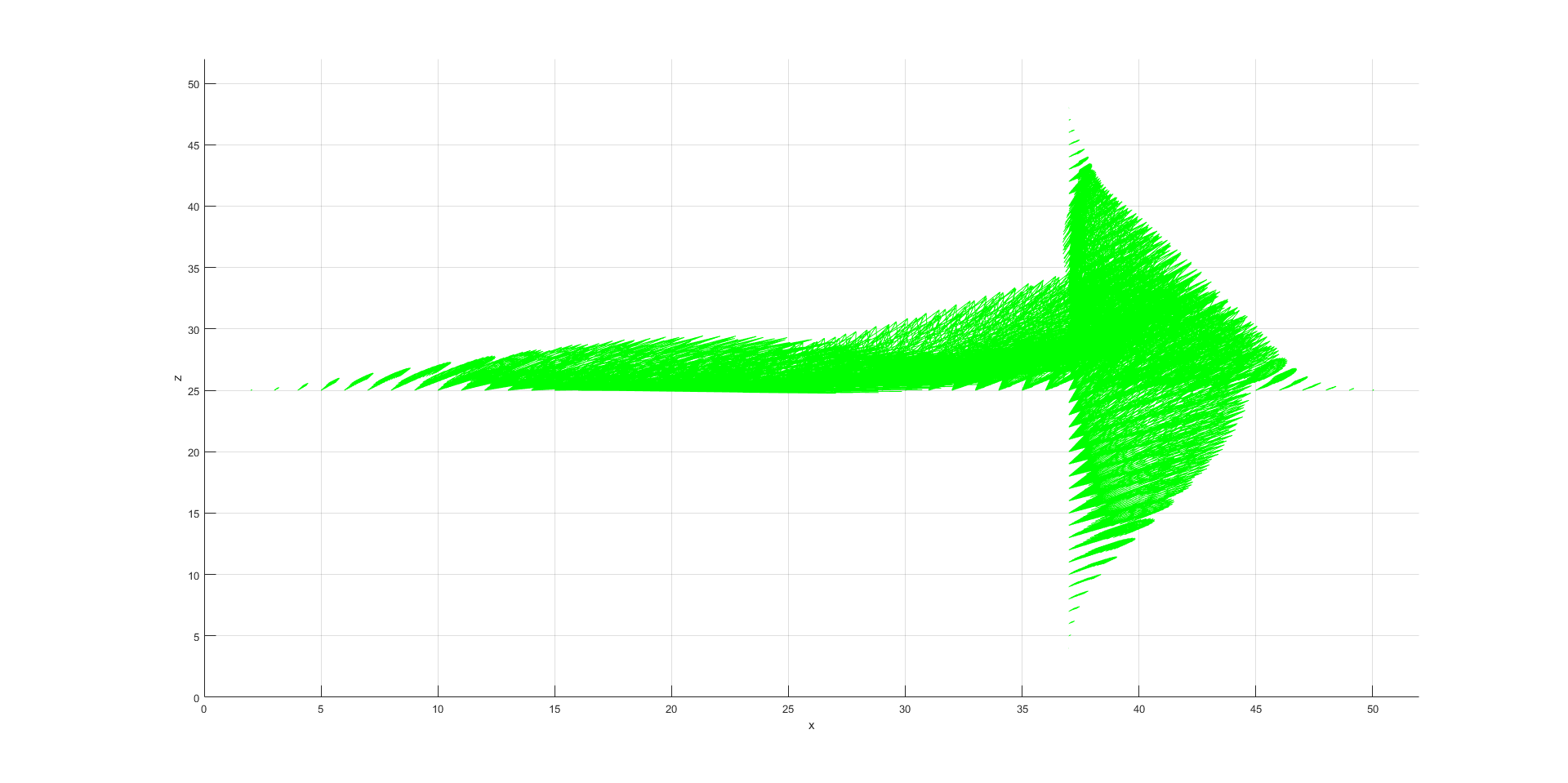}}
	\caption{Given ground truth (GT)}\label{ThreeDimGT}
\end{figure}

This example demonstrate revised VP achieves better solutions over the original VP proposed in \cite{ChenXi}. Given a 3D grid in black $\pmb{\Phi}$ and its displacement vector field $\pmb{U} = \pmb{\Phi} - \pmb{id}$ in green, in Fig.\ref{ThreeDimGT}, $\pmb{\Phi}$ is manually built by multiple times of applying cutoff rotation and displacement over the domain $\mathrm{\Omega}=[1, 51]^{3}$ such that det$\nabla\pmb{\Phi}>0.1734$, i.e., $\pmb{\Phi}$ has non-folding grids. For a cleaner visualization, all 3D grids and their displacement fields are only plotted on the 25-th frame along $z$-axis and the 37-th frame along $x$-axis. Define $f_o=$det$\nabla\pmb{\Phi}$ and $\pmb{g}_o=\nabla\times\pmb{\Phi}$ as the prescribed JD and curl.

\begin{figure}[H]	
	\subfigure[$\pmb{\Phi}_{ovp}$]{\includegraphics[width=3.4cm,height=3.4cm]{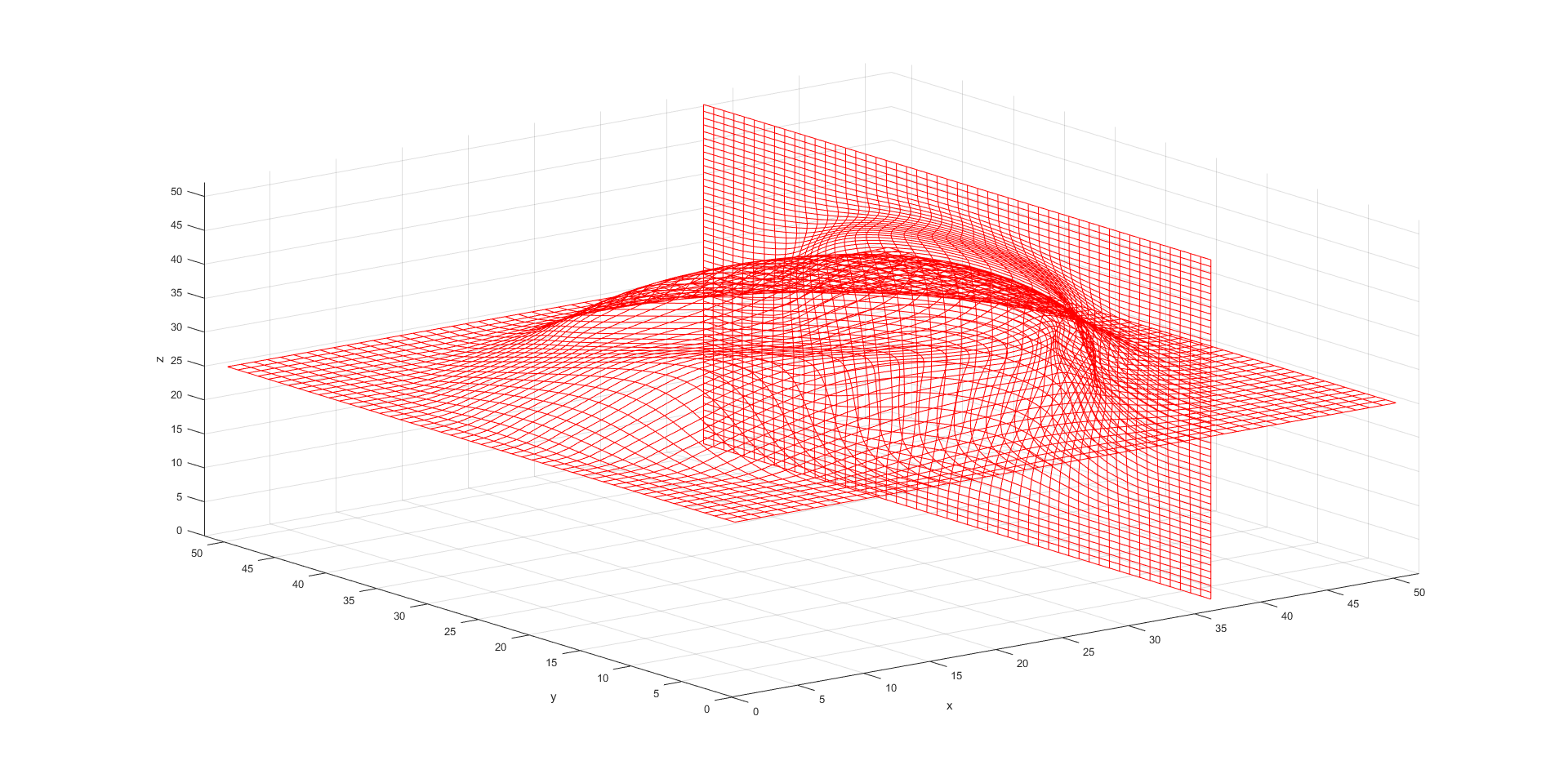}}
	\subfigure[$\pmb{\Phi}_{ovp}$ ($xz$-plain)]{\includegraphics[width=3.4cm,height=3.4cm]{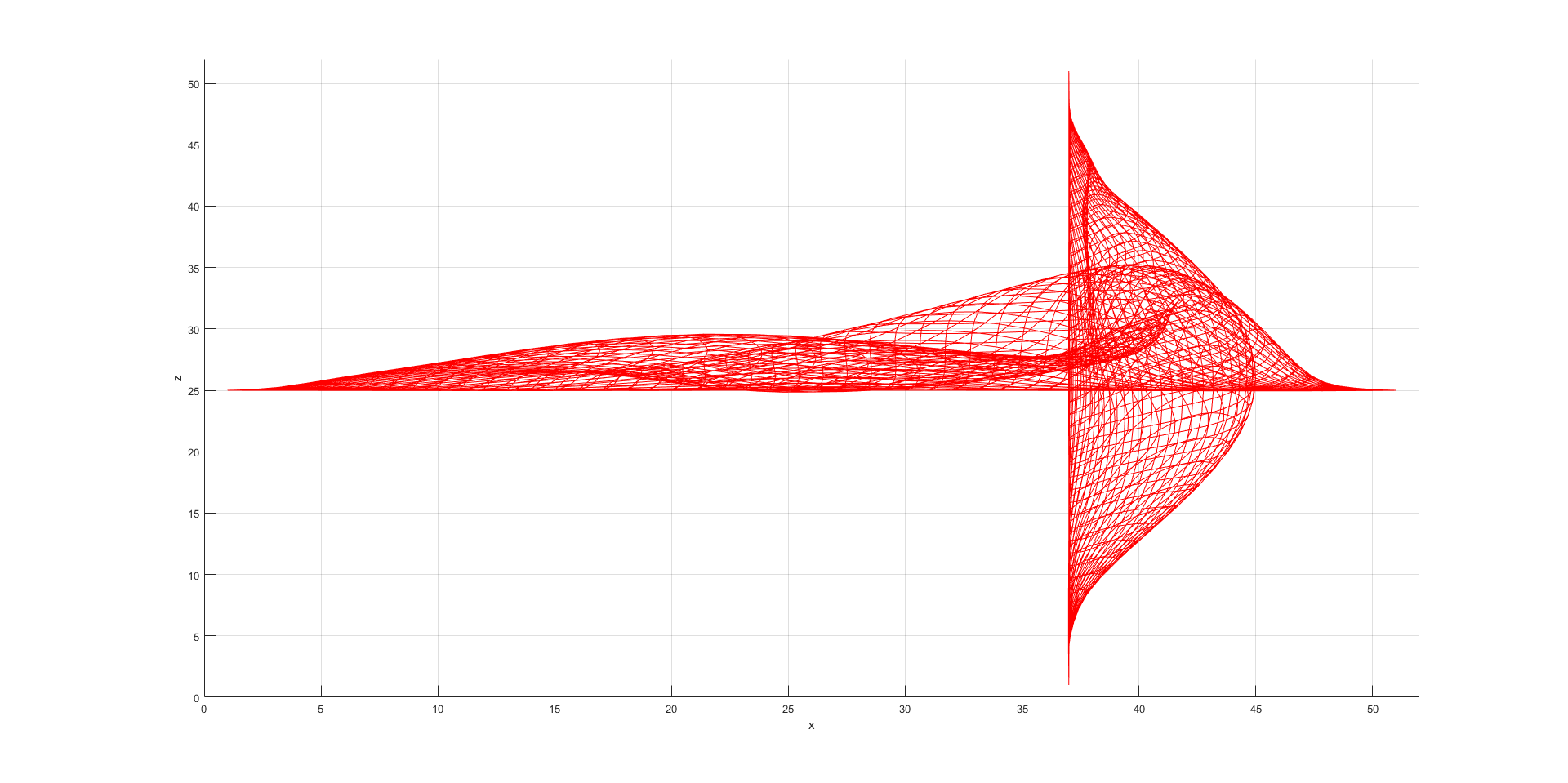}}
	\subfigure[$\pmb{U}_{ovp}$]{\includegraphics[width=3.4cm,height=3.4cm]{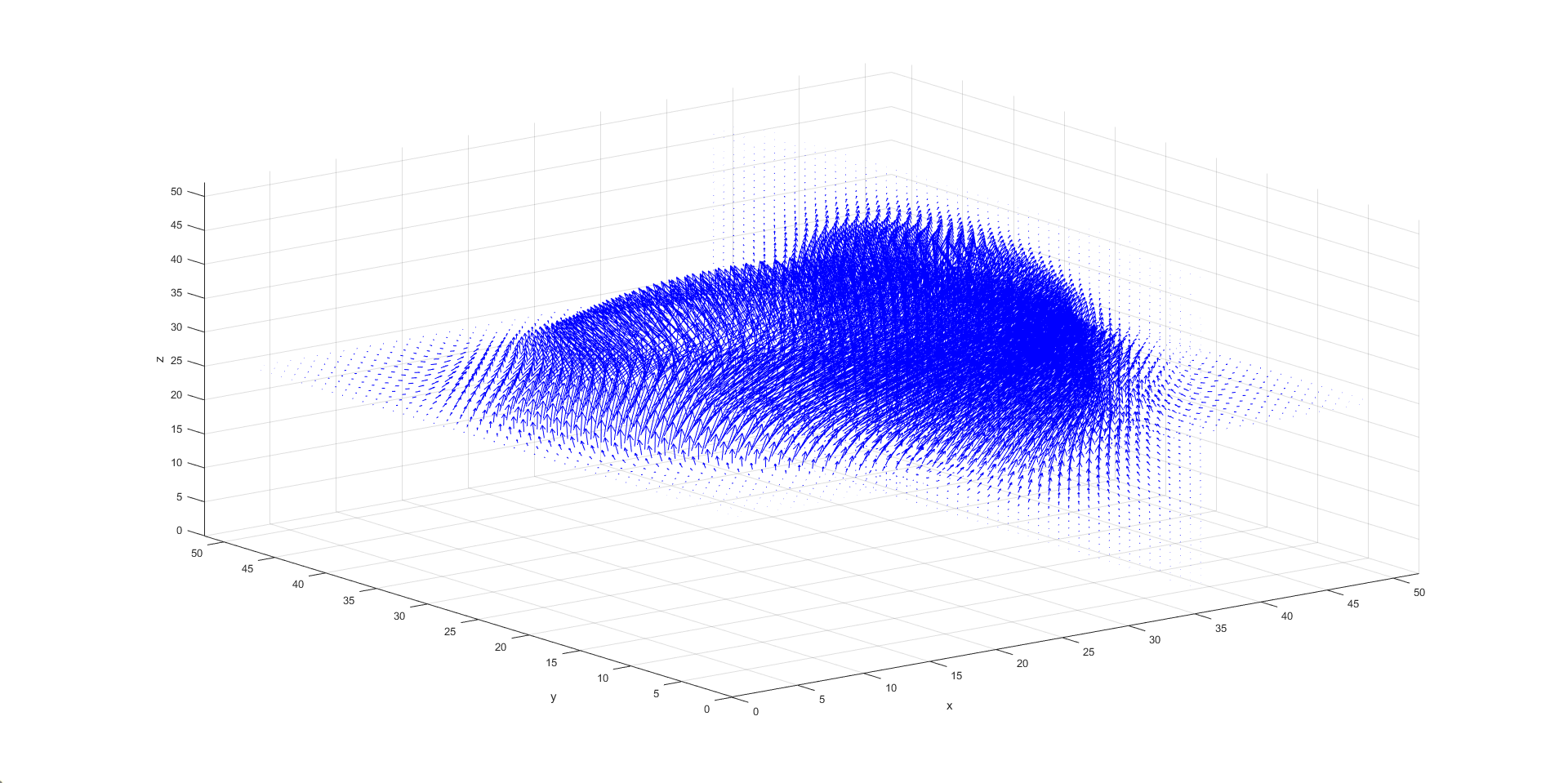}}
	\subfigure[$\pmb{U}_{ovp}$ ($xz$-plain)]{\includegraphics[width=3.4cm,height=3.4cm]{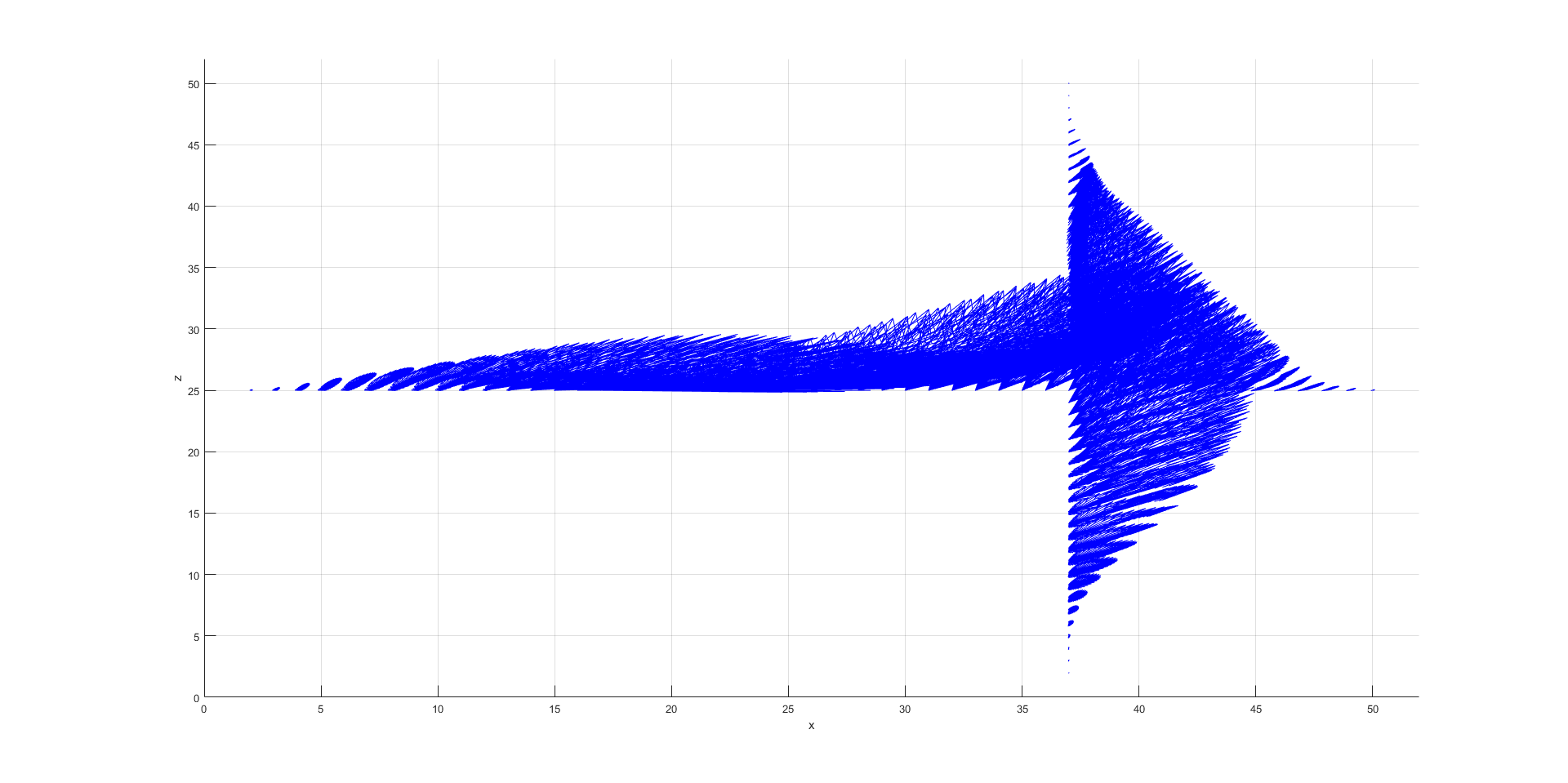}}
	
	\subfigure[$\pmb{\Phi}_{ovp}$ vs $\pmb{\Phi}$]{\includegraphics[width=3.4cm,height=3.4cm]{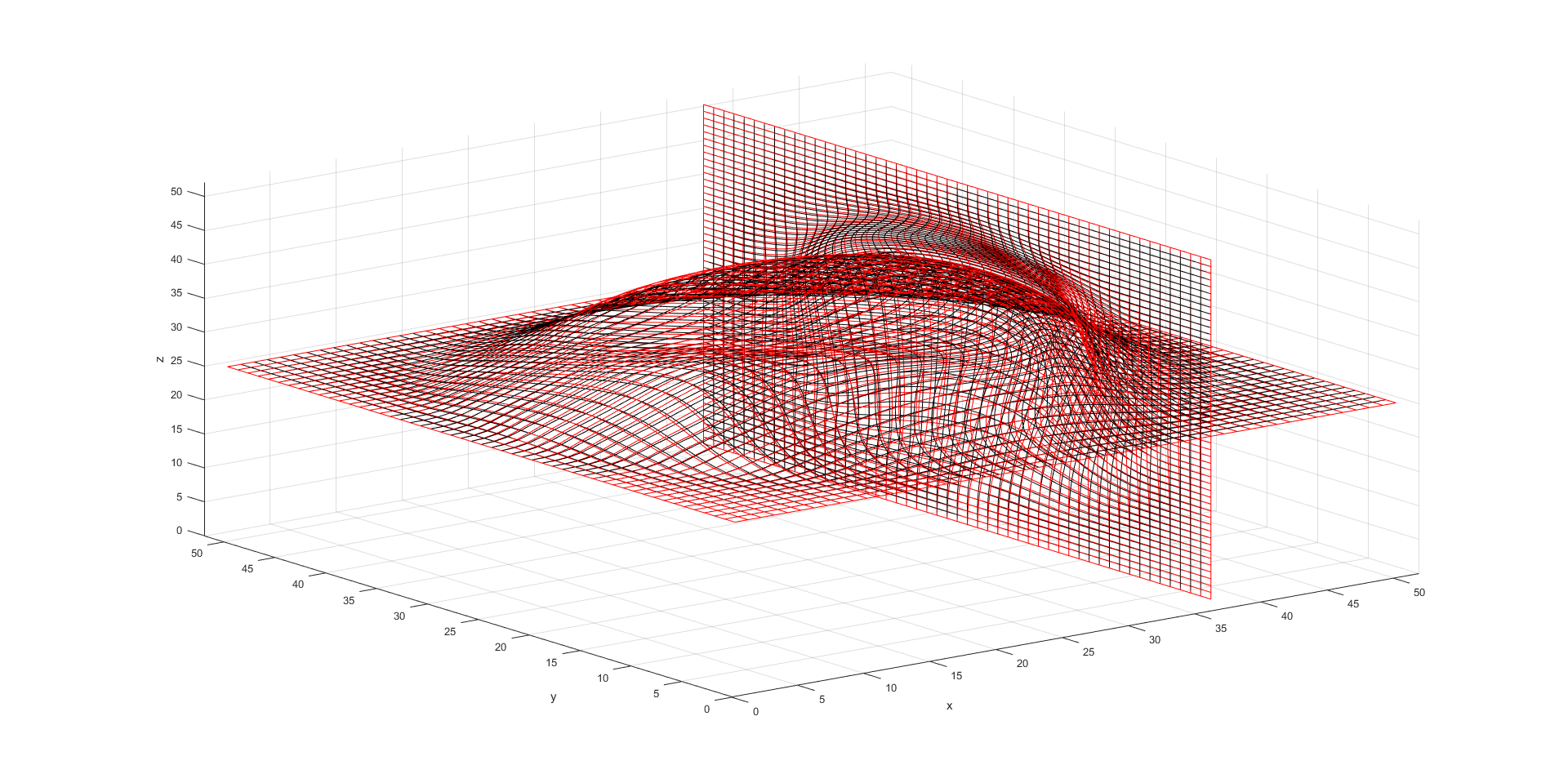}}
	\subfigure[$\pmb{\Phi}_{ovp}$ vs $\pmb{\Phi}$ ($xz$-plain)]{\includegraphics[width=3.4cm,height=3.4cm]{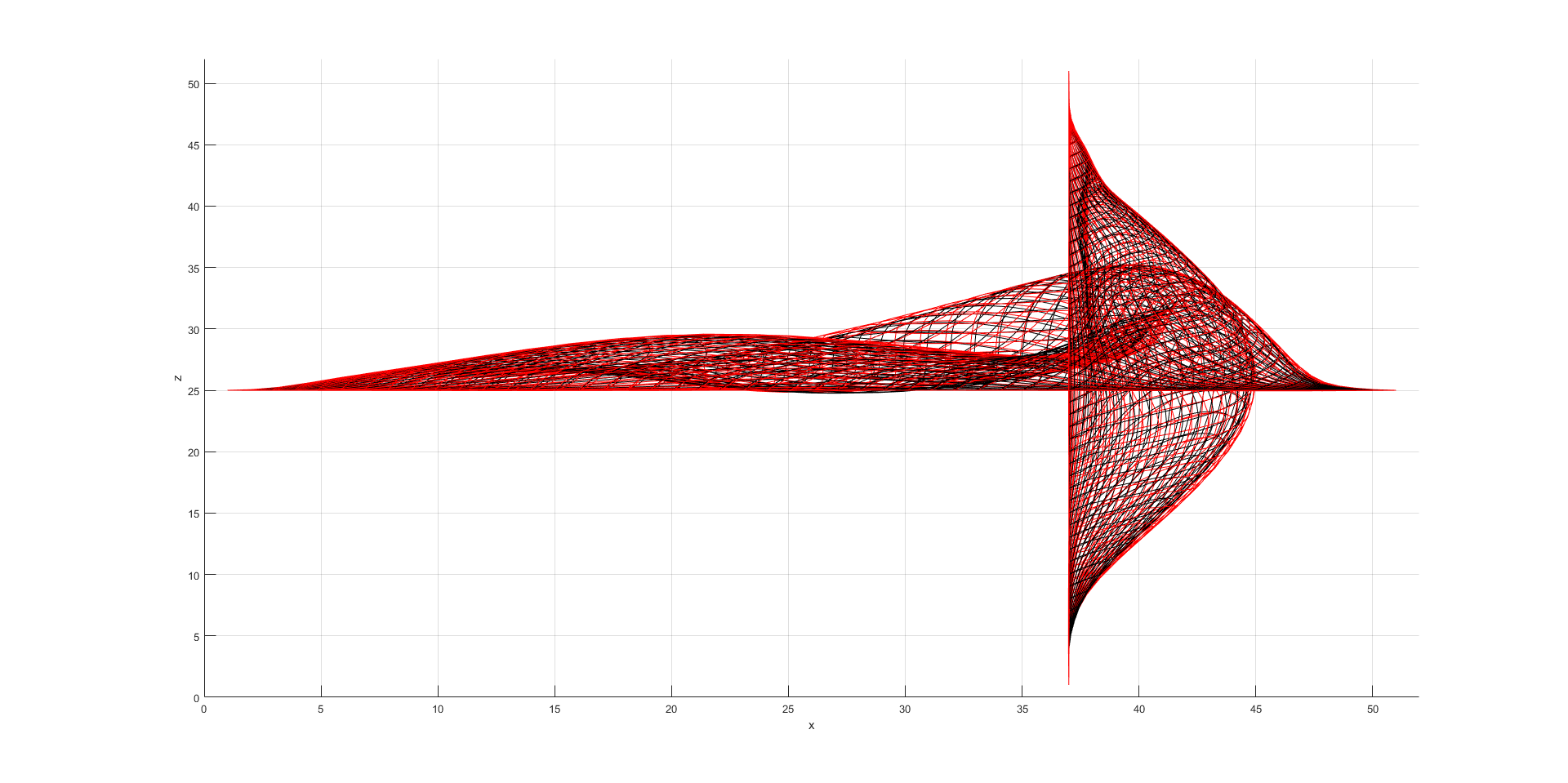}}
	\subfigure[$\pmb{U}_{ovp}$ vs $\pmb{U}$]{\includegraphics[width=3.4cm,height=3.4cm]{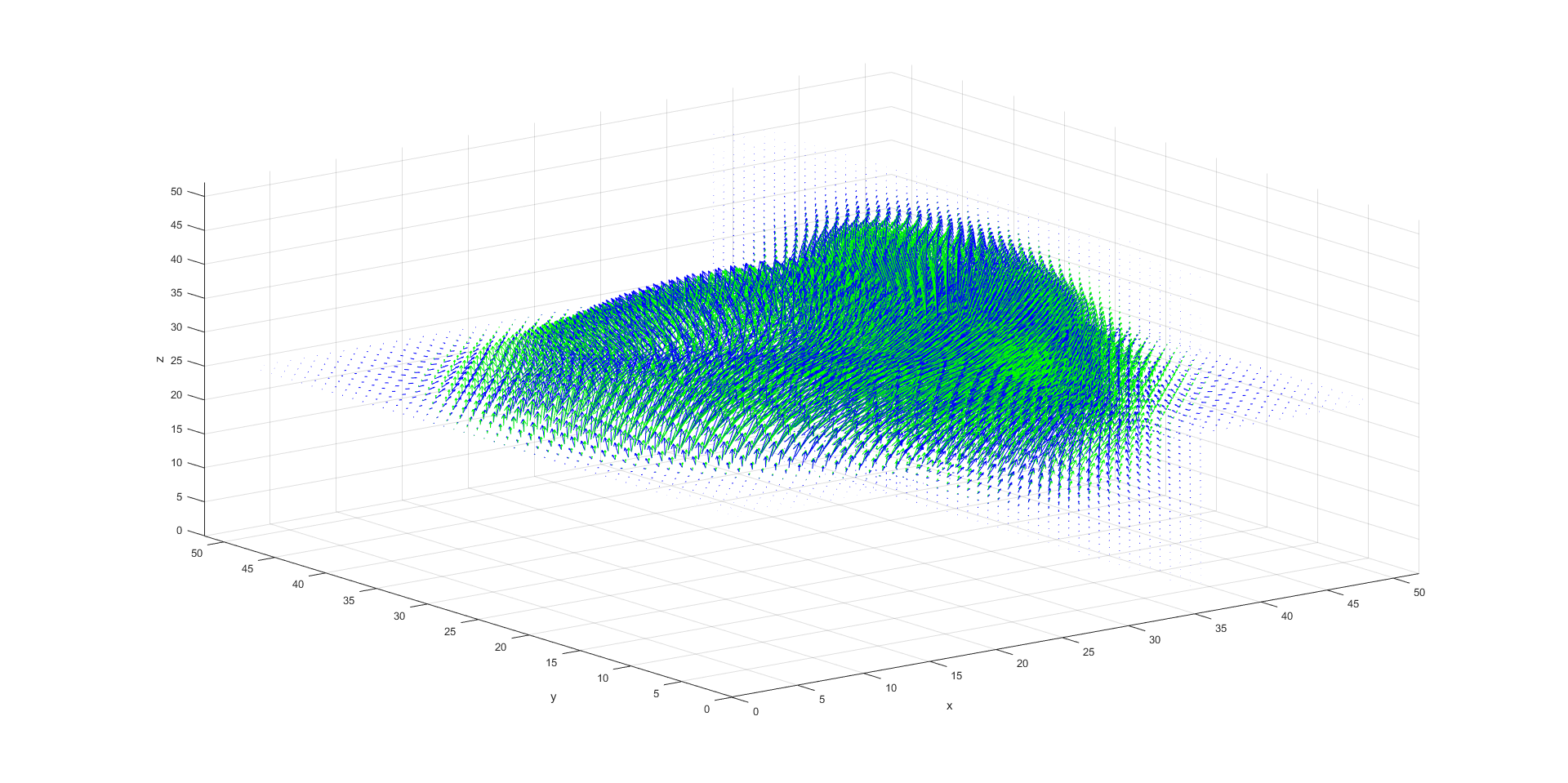}}
	\subfigure[$\pmb{U}_{ovp}$ vs $\pmb{U}$ ($xz$-plain)]{\includegraphics[width=3.4cm,height=3.4cm]{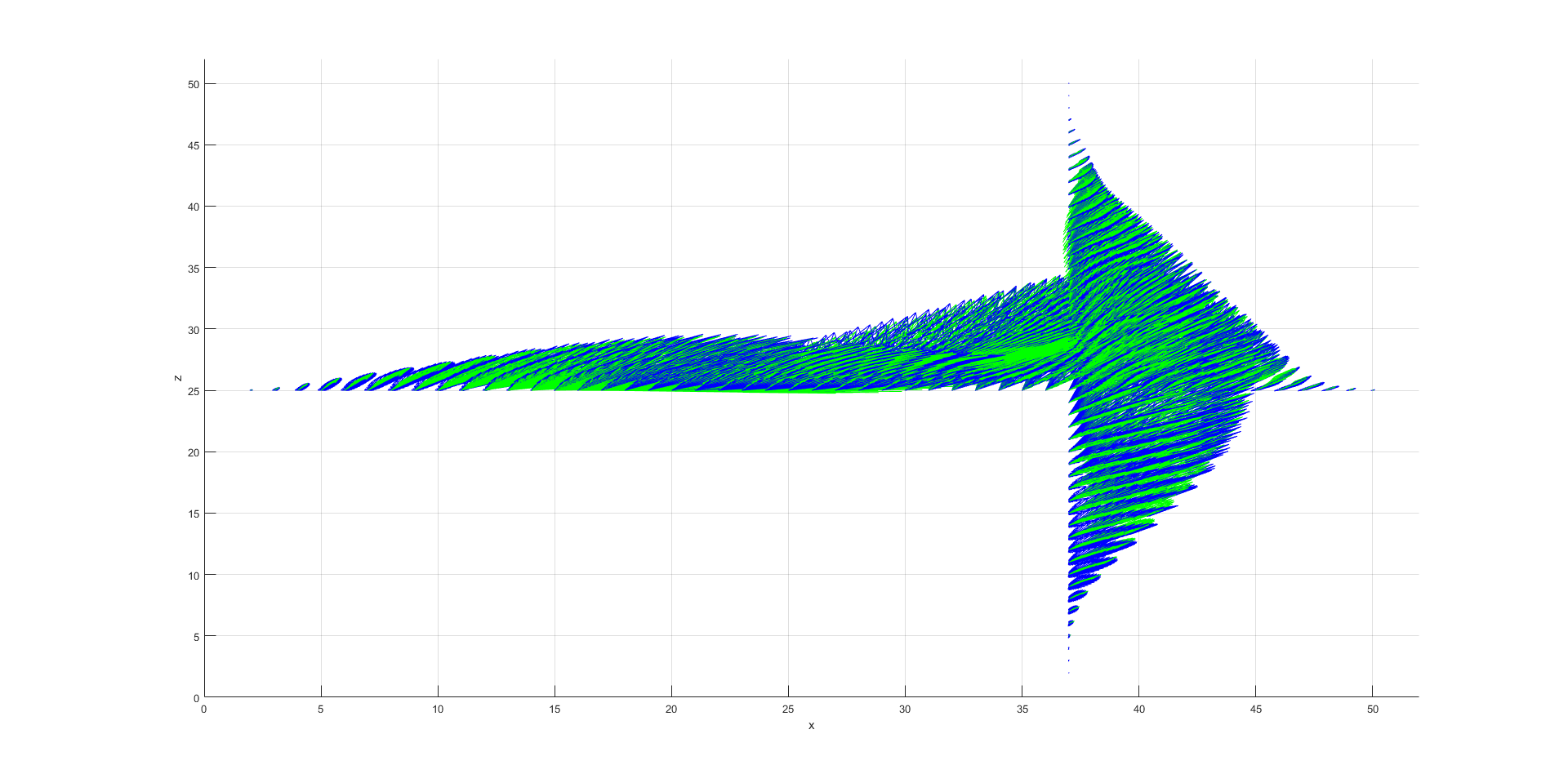}}
	\caption{Solution by original VP --- $\pmb{\Phi}_{ovp}$ compare to GT --- $\pmb{\Phi}$}\label{ThreeDimOldVPvsGT}
\end{figure}

\begin{figure}[H]
	\subfigure[$\pmb{\Phi}_{rvp}$]{\includegraphics[width=3.4cm,height=3.4cm]{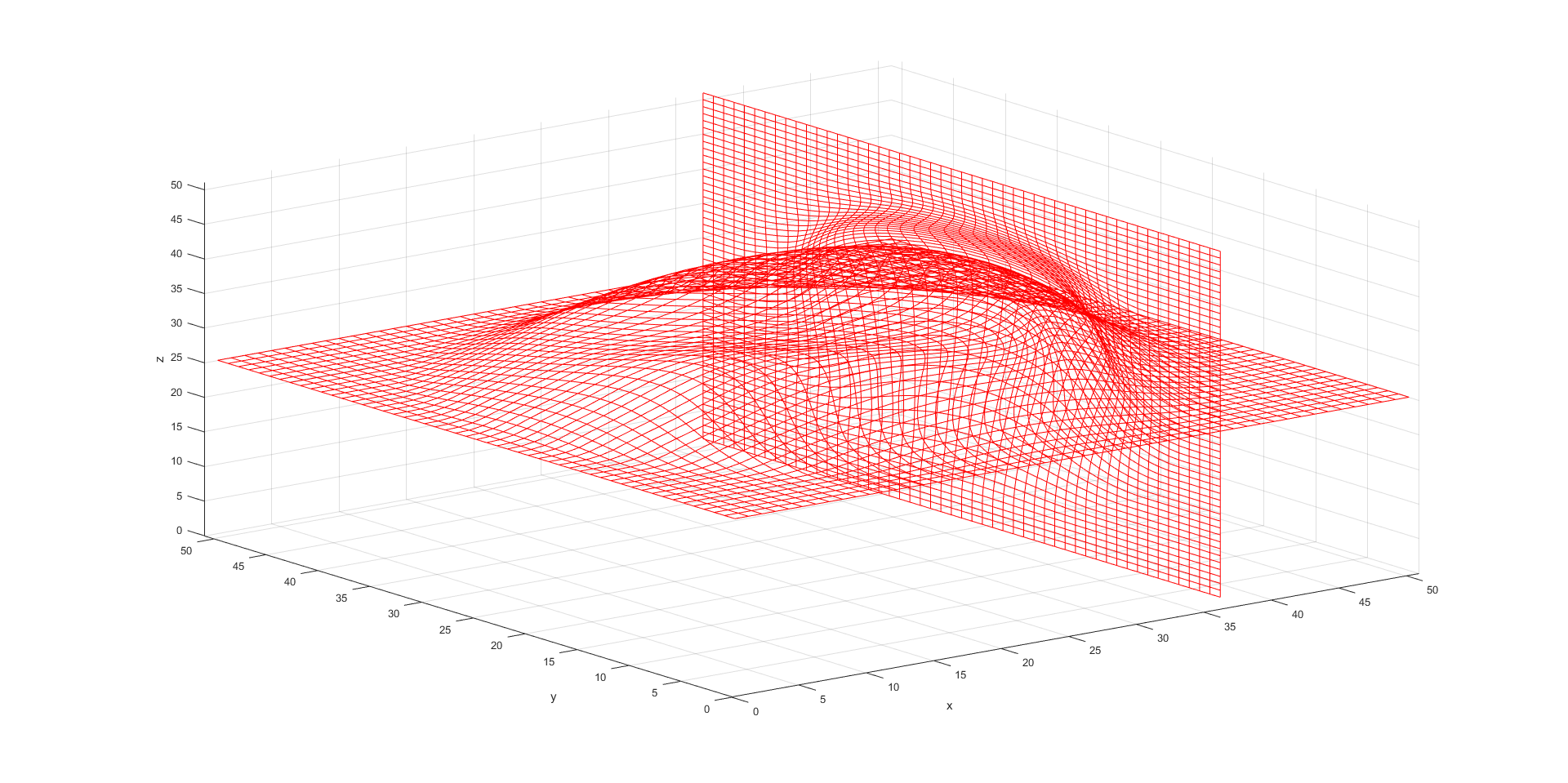}}
	\subfigure[$\pmb{\Phi}_{rvp}$ ($xz$-plain) ]{\includegraphics[width=3.4cm,height=3.4cm]{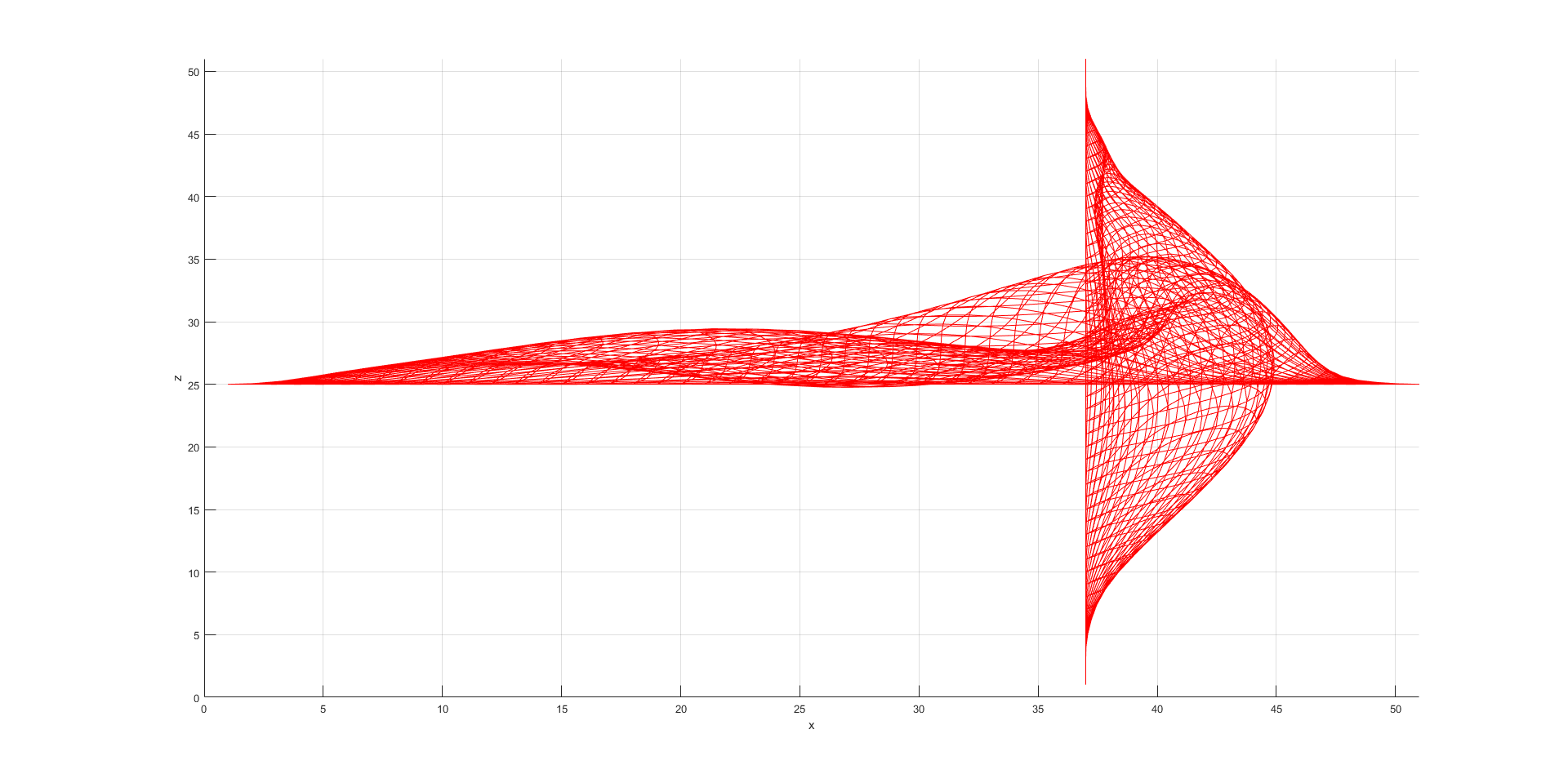}}
	\subfigure[$\pmb{U}_{rvp}$]{\includegraphics[width=3.4cm,height=3.4cm]{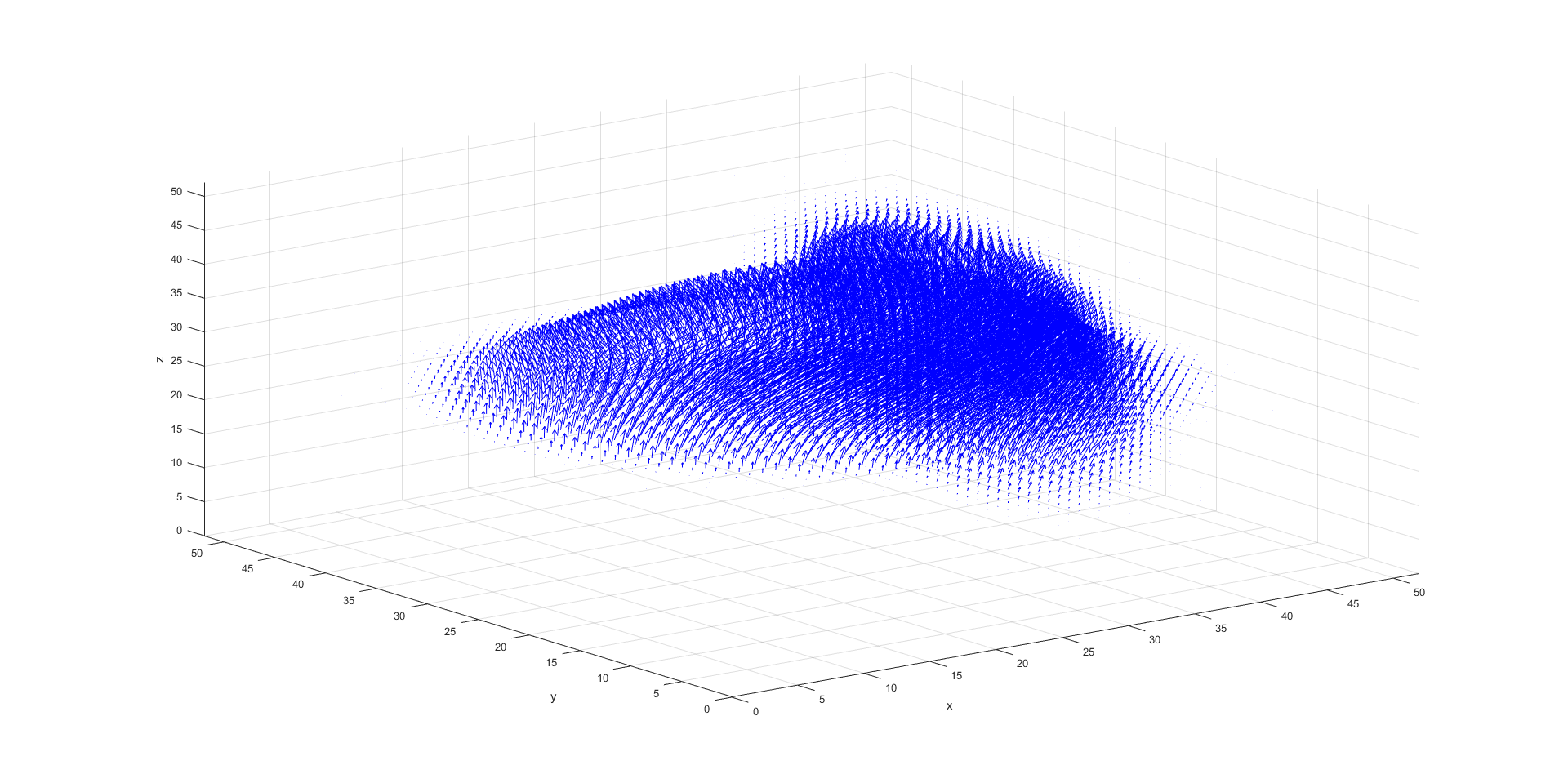}}
	\subfigure[$\pmb{U}_{rvp}$ vs $\pmb{U}$ ($xz$-plain)]{\includegraphics[width=3.4cm,height=3.4cm]{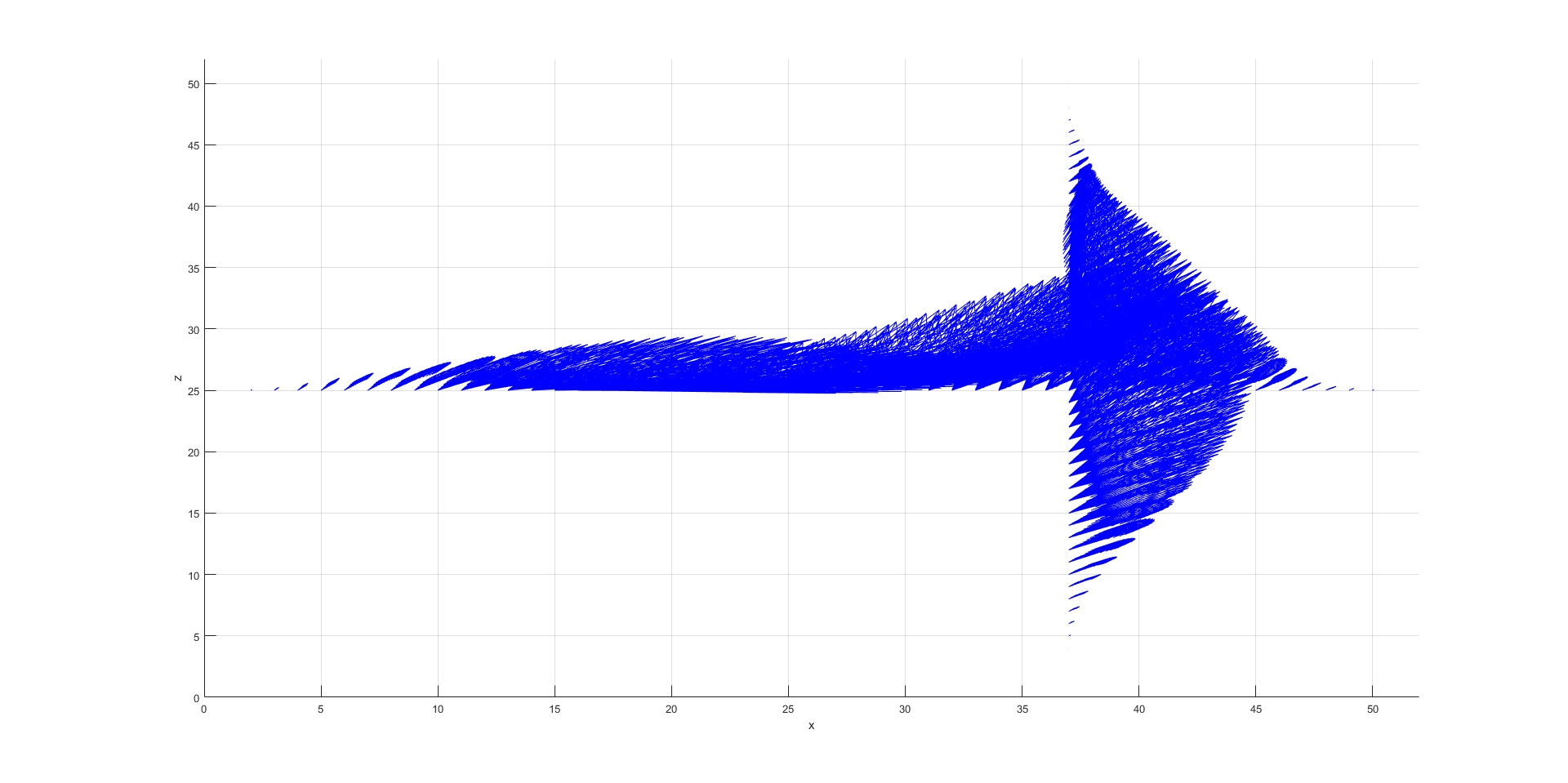}}
	
	\subfigure[$\pmb{\Phi}_{rvp}$ vs $\pmb{\Phi}$]{\includegraphics[width=3.4cm,height=3.4cm]{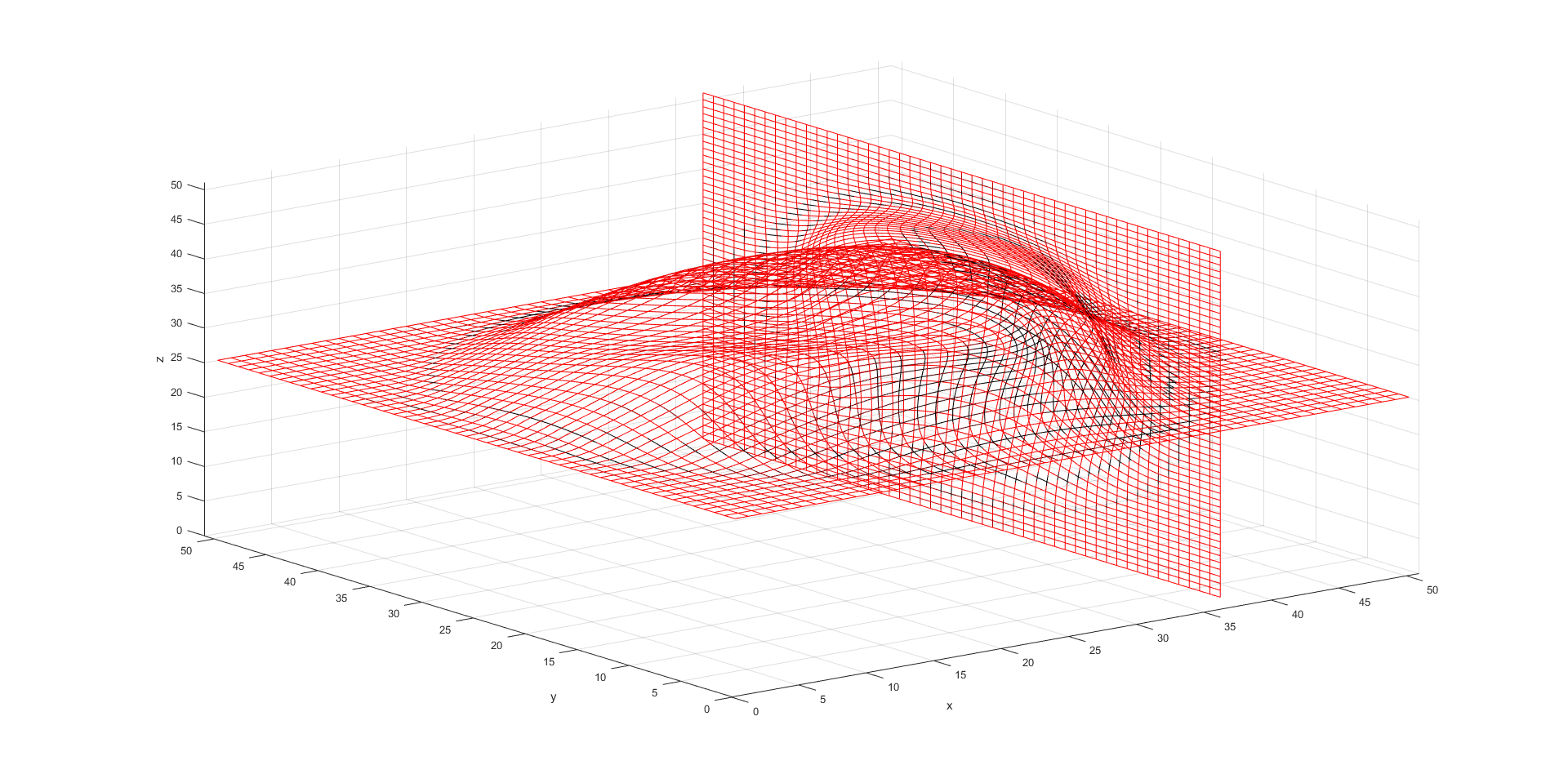}}
	\subfigure[$\pmb{\Phi}_{rvp}$ vs $\pmb{\Phi}$ ($xz$-plain)]{\includegraphics[width=3.4cm,height=3.4cm]{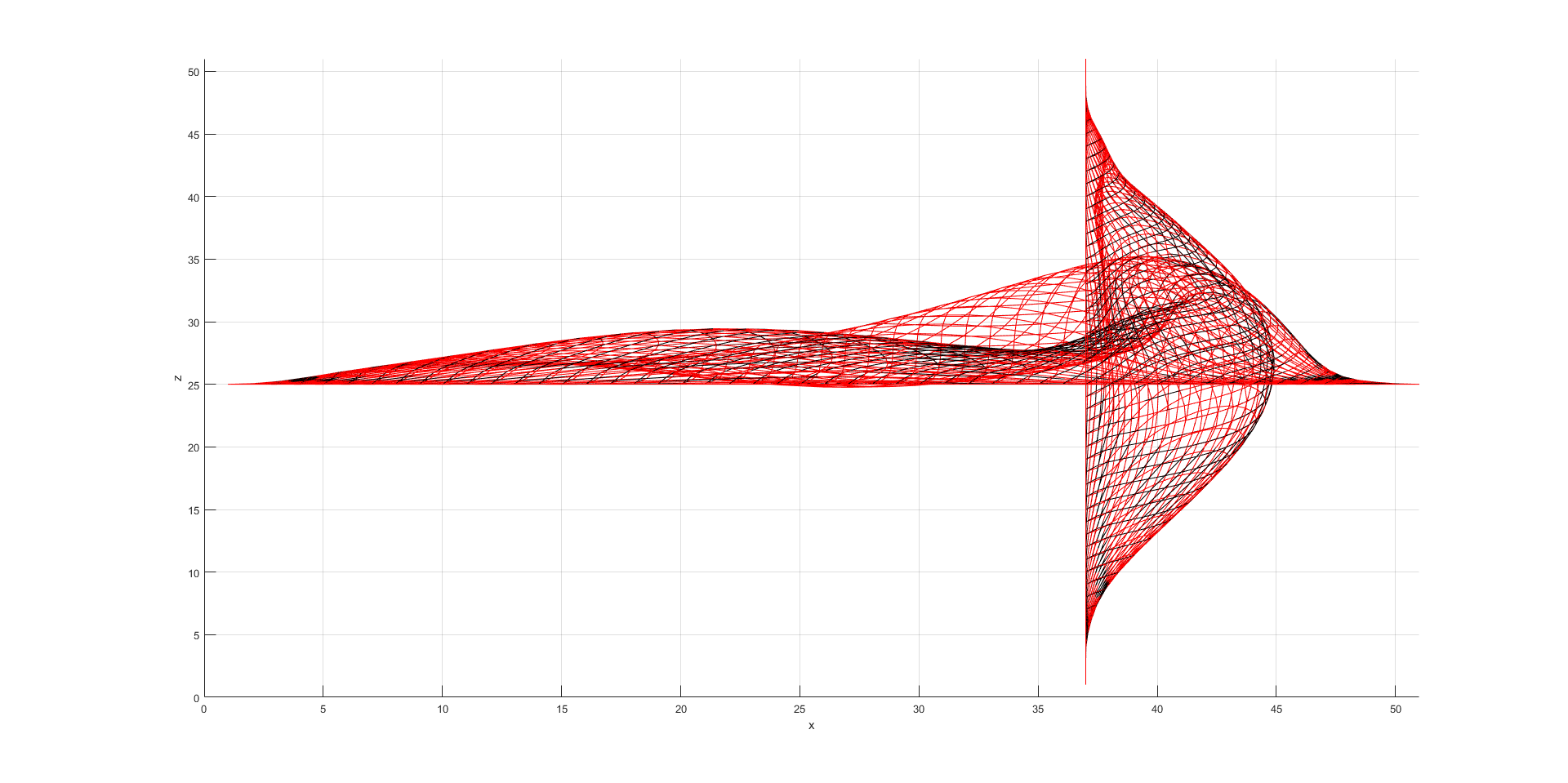}}
	\subfigure[$\pmb{U}_{rvp}$ vs $\pmb{U}$]{\includegraphics[width=3.4cm,height=3.4cm]{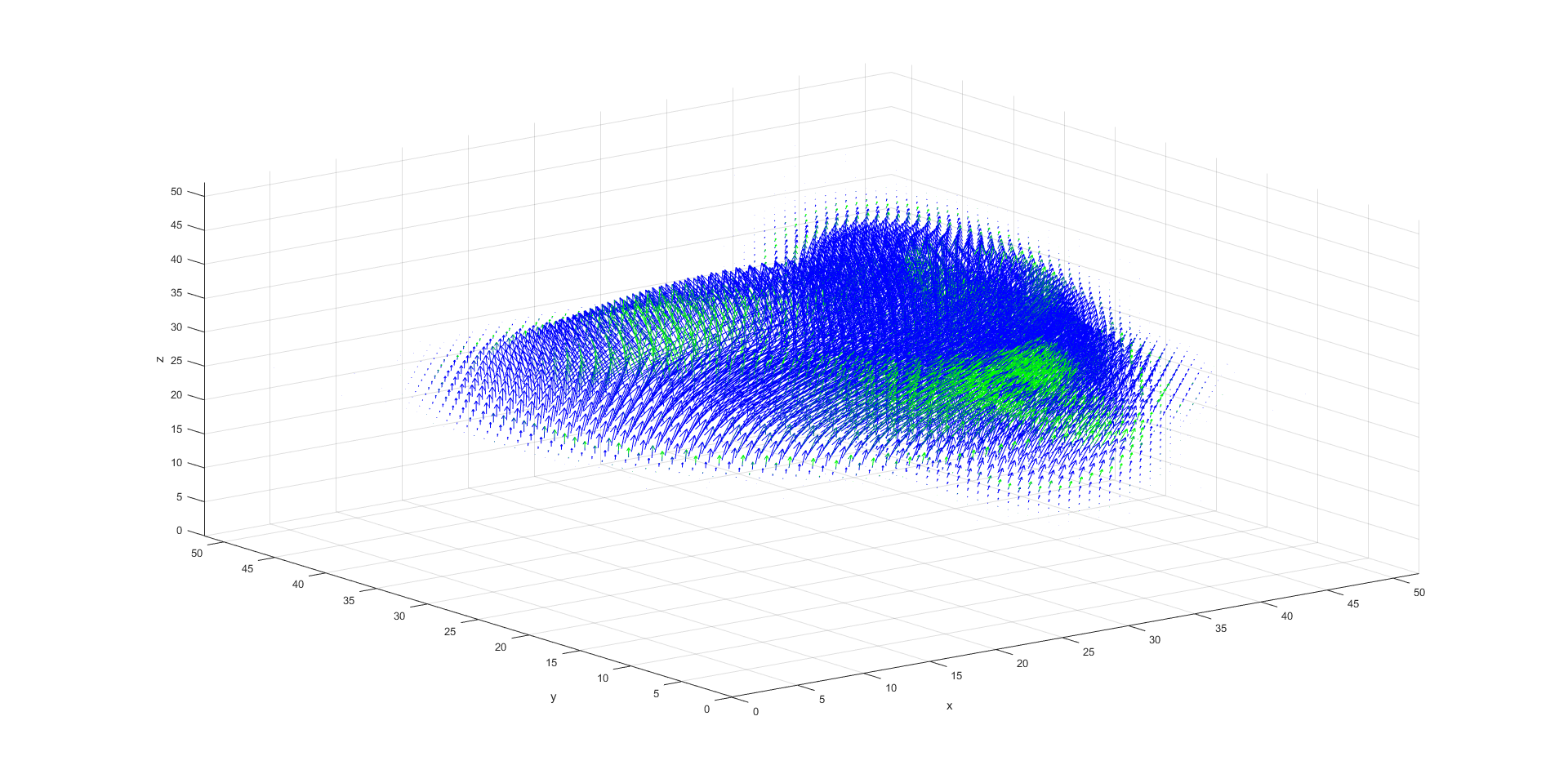}}
	\subfigure[$\pmb{U}_{rvp}$ vs $\pmb{U}$ ($xz$-plain)]{\includegraphics[width=3.4cm,height=3.4cm]{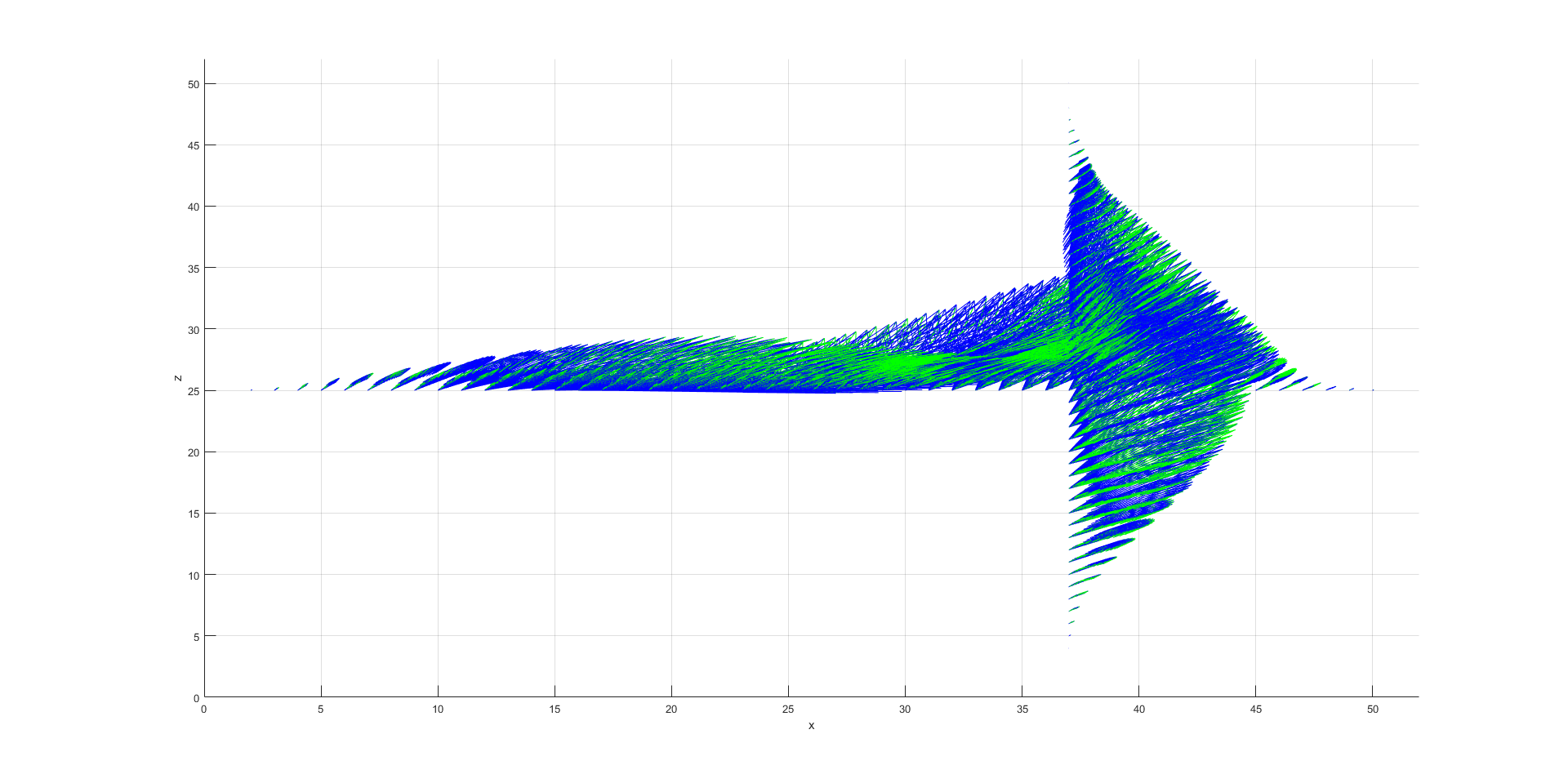}}
	\caption{Solution by revised VP --- $\pmb{\Phi}_{rvp}$ compare to GT --- $\pmb{\Phi}$}\label{ThreeDimNewVPvsGT}
\end{figure}

Fig.\ref{ThreeDimOldVPvsGT}(a-d) is the solution of original VP and is superposed on the GT in Fig.\ref{ThreeDimOldVPvsGT}(e-h). It can be visually seen that red grid lines of the solution by original VP do not line-up well with GT in black grid, in Fig.\ref{ThreeDimOldVPvsGT}(e), i.e., a better solution would have covered more black grid. The displacement vector field of original VP, in blue, is not pointing so close to the directions where the displacement vector field of GT points to, in Fig.\ref{ThreeDimOldVPvsGT}(h). Comparing to the solution of revised VP, in Fig.\ref{ThreeDimNewVPvsGT}(a-d), that also is superposed on GT, in Fig.\ref{ThreeDimNewVPvsGT}(e-h), and there is only a small portion of the black grid and the green vectors stay uncovered by the solution of revised VP, which means revised VP provides much better solutions over the original VP. This observation is confirmed by measurements in the following Table.\ref{tbl34}. It also recorded the revised VP had reached the $ratio$-tolerance, $10^{-4}$\%, within lesser iterations and computational time. 


%
%
%
%
%
%
%
%
%
%
%

\begin{table}[h]
	\small{
		\begin{center}
			\begin{tabular}{|c|c|c|ccc|}
				\hline	
				\multirow{2}{*}{Solution} & \multirow{2}{*}{sec $|$ iteration} & \multirow{2}{*}{$ratio$}  &\multicolumn{3}{c|}{max differences of Fig.\ref{ThreeDimOldVPvsGT}(e) \& \ref{ThreeDimNewVPvsGT}(e)} \\  \cline{4-6}
				&      &       &$|$det$\nabla(\_)-f_{o}|$ & $||\nabla\times(\_)-\pmb{g}_{o}||_{2}$ &  $|| (\_)-\pmb{\Phi}||_{2}$ \\ \hline
				{\tt $\pmb{\Phi}_{ovp}$} &1765.89 $|$ 3016 & $0.0052\%$  & 0.7279  &0.0866 &  $1.9454$ \\
				{\tt$\pmb{\Phi}_{rvp}$} & 788.21 $|$ 576& $0.0009\%$  & 0.0104 & 0.0136 & $0.0403$ \\ 
				\hline
			\end{tabular}
			\caption{Performance of Fig.\ref{ThreeDimOldVPvsGT} \& \ref{ThreeDimNewVPvsGT}}\label{tbl34}
		\end{center}
	}
\end{table}

This example summarizes a three-fold improvement of the revised VP: (1) the revised VP is effective and capable in generating 3D grids with prescribed JD and curl; (2) the revised VP finds more accurate solutions in terms of the prescription of JD and curl (it is more obvious in 3D examples) compare to the original VP; (3) the revised VP reaches the desired tolerance faster than the original VP (in fact, the computational cost of the revised VP for an effective iteration are very close the original VP, where main extra computations occur on the interpolation of step 9 in the provided algorithm). 

%
%
%
%
%
%
%
%
%
%
%

\subsection{$\mathbf{Example}$: Reconstruction of  Inverse}\label{eg1}
\begin{figure}[H]
	\centering
	\subfigure[$\pmb{B}_{o}$]{\includegraphics[width=3.4cm,height=3.4cm]{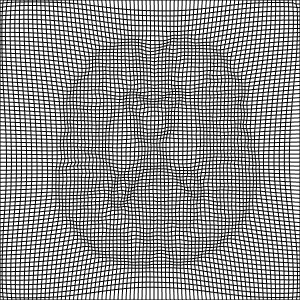}}
	\subfigure[$\pmb{u}$]{\includegraphics[width=3.4cm,height=3.4cm]{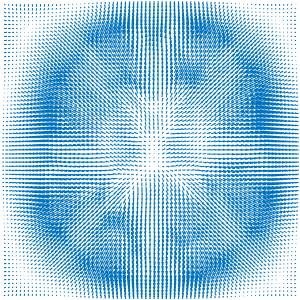}}
	\subfigure[$\pmb{D}$]{\includegraphics[width=3.4cm,height=3.4cm]{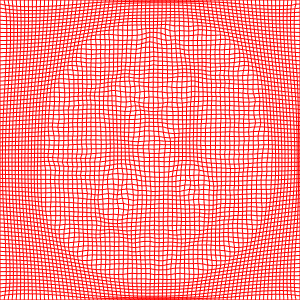}}
	\subfigure[$\pmb{D}\circ\pmb{B}_{o}$]{\includegraphics[width=3.4cm,height=3.4cm]{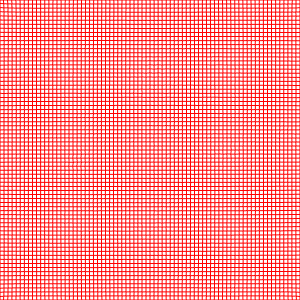}}
	
	\subfigure[$\pmb{D}\circ\pmb{B}_{o}$ vs $\pmb{id}$]{\includegraphics[width=3.4cm,height=3.4cm]{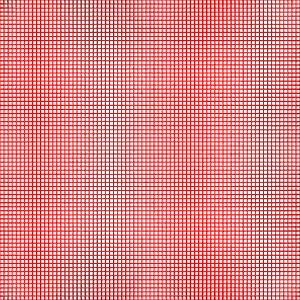}}
	\subfigure[$\pmb{B}_{o}\circ\pmb{D}$ vs $\pmb{id}$]{\includegraphics[width=3.4cm,height=3.4cm]{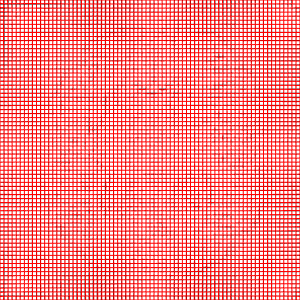}}	
	\caption{Reconstruct inverse under left-translation}\label{debrain}
\end{figure}

Given a brain-like grid, $\pmb{B}_{o}$, Fig.\ref{debrain}(a). Because det$\nabla\pmb{id}=1$ and $\nabla\times\pmb{id}=\pmb{0}$ for the identity map $\pmb{id}$, so $1=f_o$, $\pmb{0}=\pmb{g}_{o}$ and $\pmb{B}_{o}=\pmb{\phi}_{o}$ are fed to the algorithm. It is expected to find a transformation that composites $\pmb{B}_{o}$ and outputs $\pmb{id}$. As Fig.\ref{debrain}(c) shows, $\pmb{D}$ is found and it left-translates $\pmb{B}_{o}$ to Fig.\ref{debrain}(d). Compared it with $\pmb{id}$ in Fig.\ref{debrain}(e), the brain-like grid $\pmb{B}_{o}$ is ``de-brained". Fig.\ref{debrain}(f) is a guessed identity by reversing the order of $\pmb{B}_{o}$ and $\pmb{D}$ in Fig.\ref{debrain}(e). This indicates that $\pmb{D}$ approximates quite accurately to $(\pmb{B}_{o})^{-1}$. $ratio$ had reduced to 0.3\%.

\subsection{$\mathbf{Example}$: Inverse Consistency}\label{eg2}
Given a grid, $\pmb{\Psi}$, of the character  ``\hyperref{https://en.wikipedia.org/wiki/Li_(surname_%E6%9D%8E)#/media/File:Li_(%E6%9D%8E).svg}{link}{Li}{Li}''
 that is used to symbolize Lie groups in Chinese, and a grid, $\pmb{\Phi}$,  of a Mona Lisa image in Fig.\ref{invconsis}. Feed det$\nabla\pmb{\Phi}=f_o$, $\nabla\times\pmb{\Phi}=\pmb{g}_{o}$ and $\pmb{\Psi}=\pmb{\phi}_{o}$ to the algorithm, then $\pmb{\phi}^{k}=\pmb{D}_{1}\circ\pmb{\Psi}$ and $\pmb{\phi}^{k}_{\pmb{m}}=\pmb{D}_{1}$ are outputted. By reversing the order, $\pmb{D}_{2}$ and $\pmb{D}_{2}\circ\pmb{\Psi}$ are found. It is expected that $\pmb{D}_{1}$ and $\pmb{D}_{2}$ to be the inverse transformation of each other. It is shown in Fig.\ref{invconsis}(k,l) where the compositions of $\pmb{D}_{1}$ and $\pmb{D}_{2}$ are close to $\pmb{id}$. 
\begin{figure}[H]
	\centering
	\subfigure[$\pmb{\Phi}$]{\includegraphics[width=3.4cm,height=3.4cm]{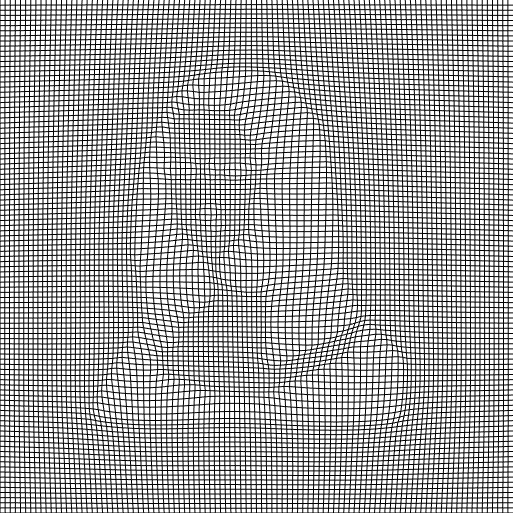}}
	\subfigure[$\pmb{\Psi}$]{\includegraphics[width=3.4cm,height=3.4cm]{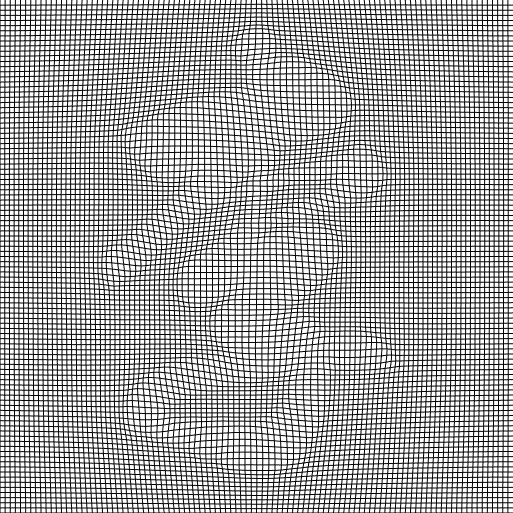}}

	\subfigure[$\pmb{u}_{1}$]{\includegraphics[width=3.4cm,height=3.4cm]{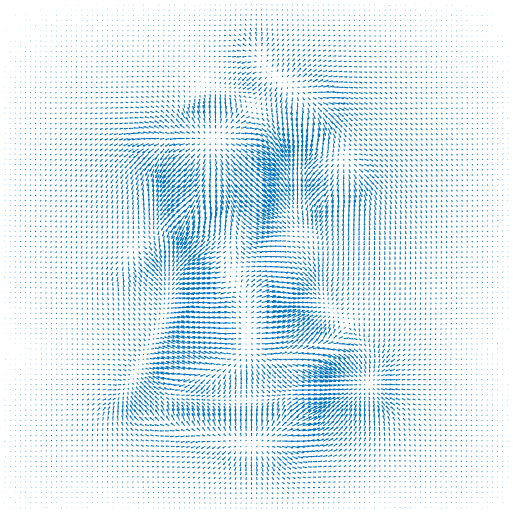}}
	\subfigure[$\pmb{D}_{1}$ ]{\includegraphics[width=3.4cm,height=3.4cm]{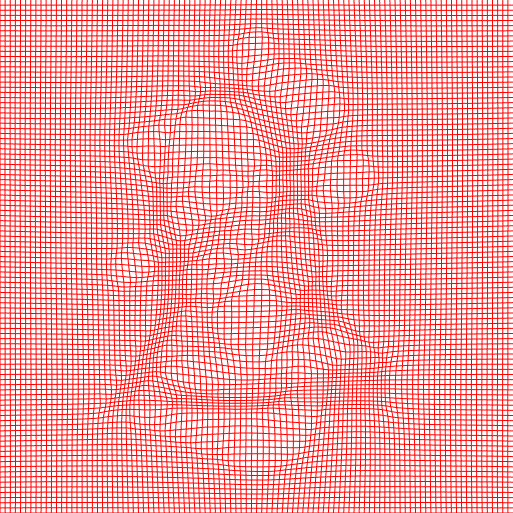}}
	\subfigure[$\pmb{D}_{1}\circ\pmb{\Phi}$ ]{\includegraphics[width=3.4cm,height=3.4cm]{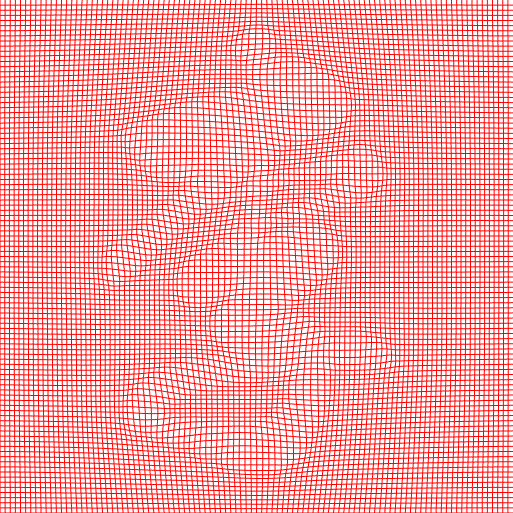}}
	\subfigure[$\pmb{D}_{1}\circ\pmb{\Phi}$ vs $\pmb{\Psi}$ ]{\includegraphics[width=3.4cm,height=3.4cm]{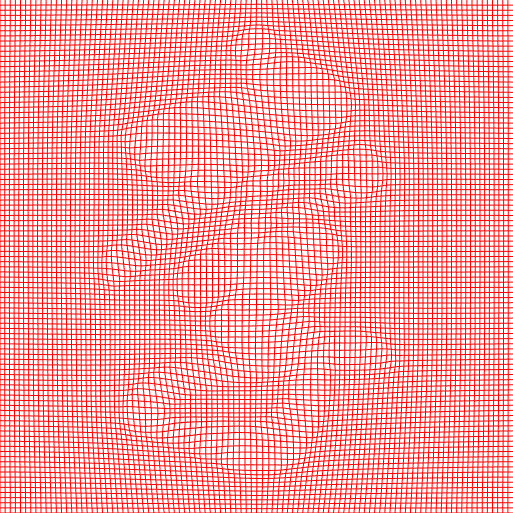}}

	\subfigure[$\pmb{u}_{2}$]{\includegraphics[width=3.4cm,height=3.4cm]{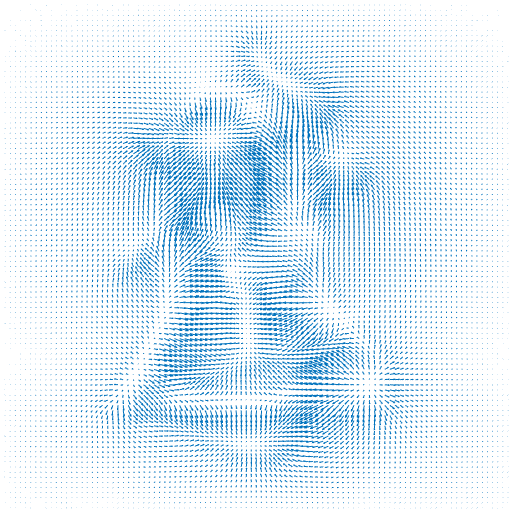}}
	\subfigure[$\pmb{D}_{2}$ ]{\includegraphics[width=3.4cm,height=3.4cm]{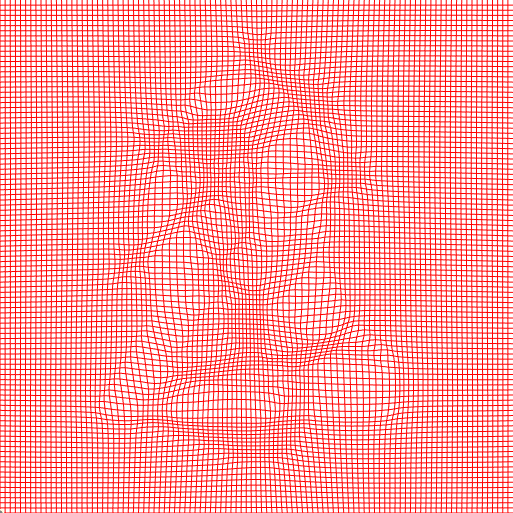}}
	\subfigure[$\pmb{D}_{2}\circ\pmb{\Psi}$ ]{\includegraphics[width=3.4cm,height=3.4cm]{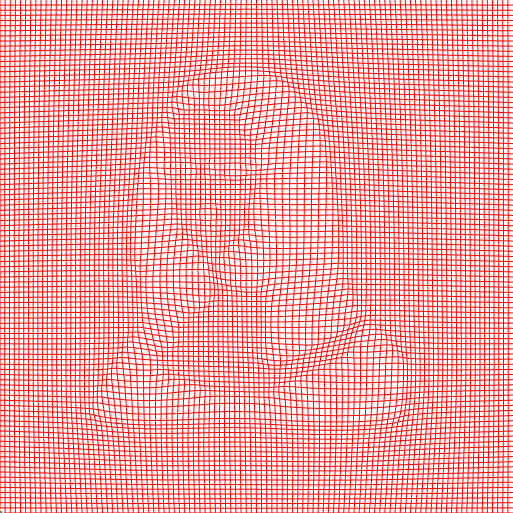}}
	\subfigure[$\pmb{D}_{2}\circ\pmb{\Psi}$ vs $\pmb{\Phi}$]{\includegraphics[width=3.4cm,height=3.4cm]{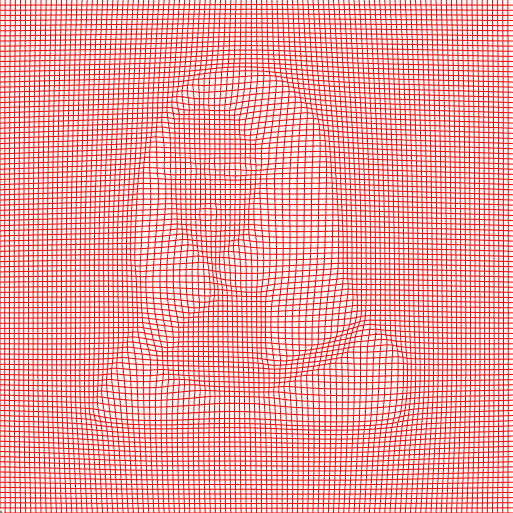}}
\end{figure}

\begin{figure}[H]
\centering	
	\subfigure[$\pmb{D}_{2}\circ\pmb{D}_{1}$ vs $\pmb{id}$ ]{\includegraphics[width=3.4cm,height=3.4cm]{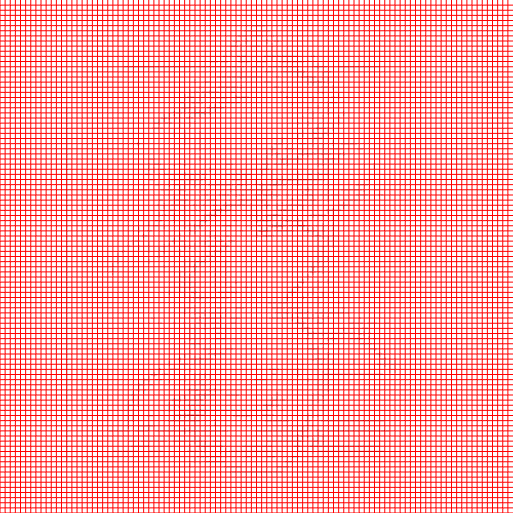}}
	\subfigure[$\pmb{D}_{1}\circ\pmb{D}_{2}$ vs $\pmb{id}$ ]{\includegraphics[width=3.4cm,height=3.4cm]{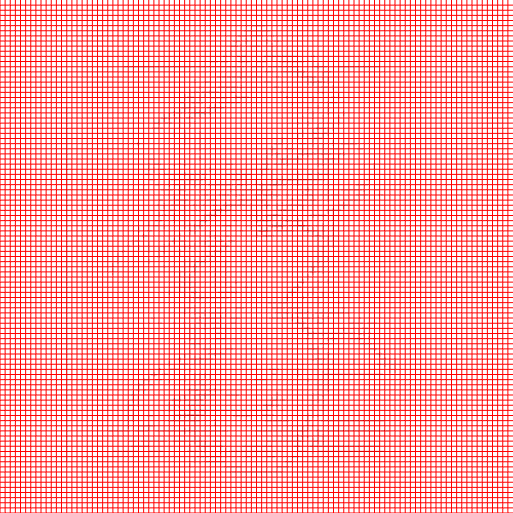}}
	\caption{Inverse consistency}\label{invconsis}
\end{figure}

\begin{table}[h]
	\small{
		\begin{center}
			\begin{tabular}{|c|c|c|c|ccc|}
				\hline
				\multirow{2}{*}{orientation} & \multirow{2}{*}{$\mathrm{\Omega}$}& \multirow{2}{*}{$ratio$} & \multirow{2}{*}{sec $|$ iteration} &\multicolumn{3}{c|}{max differences of  Fig.\ref{invconsis}(f, j)} \\  \cline{5-7} 
				&      &    &    &$|$det$\nabla(\_)|$&$|\nabla\times(\_)|$&  $|| (\_) ||_{2}/|\mathrm{\Omega}|$   \\ \hline
				{\tt$\pmb{\Phi}$ to $\pmb{\Psi}$} &$[1,101]^2$ & $0.009\%$  & 172.81 $|$ 16623  & 0.0160 &0.0211& $2.1139*10^{-6}$ \\ 
				{\tt $\pmb{\Psi}$ to $\pmb{\Phi}$} &$[1,101]^2$ & $0.009\%$ & 168.65 $|$ 16624  & 0.0174  &0.0125 &  $2.1069*10^{-6}$\\
				\hline
			\end{tabular}
			\caption{Performance of Fig.\ref{invconsis}}\label{tbl3}
		\end{center}
	}
\end{table}

\subsection{$\mathbf{Example}$: Transitivity}\label{eg3}
A problem of transitivity is presented in this example. $\pmb{C}$,  $\pmb{P}$ and $\pmb{R}$ are grids of circle, pentagon and rectangle shapes, respectively, shown in Fig.\ref{trans}(a, b, c). Firstly, $\pmb{D}_{1}$ from $\pmb{C}$ to $\pmb{P}$, $\pmb{D}_{2}$ from $\pmb{P}$ to $\pmb{R}$ and $\pmb{D}_{3}$ from $\pmb{C}$ to $\pmb{R}$ are found and shown in Fig.\ref{trans}(e, h, m), respectively; second, form $\pmb{D}$ Fig.\ref{trans}(j) by left-translating $\pmb{D}_{1}$ with $\pmb{D}_{2}$, then compare it with $\pmb{D}_{3}$ as shown in Fig.\ref{trans}(o). Fig.\ref{trans} has shown what is expected in the sense of transitivity. 
\begin{figure}[H]
	\centering
	\subfigure[$\pmb{C}$]{\includegraphics[width=3.4cm,height=3.4cm]{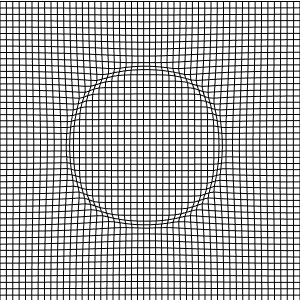}}
	\subfigure[$\pmb{P}$]{\includegraphics[width=3.4cm,height=3.4cm]{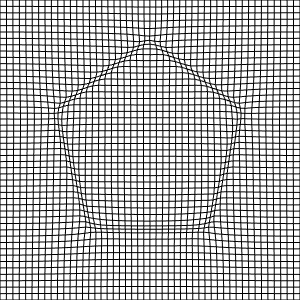}}
	\subfigure[$\pmb{R}$]{\includegraphics[width=3.4cm,height=3.4cm]{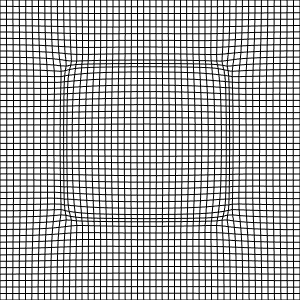}}
\end{figure}

\begin{figure}[H]
\centering	
	\subfigure[$\pmb{u}_{1}$]{\includegraphics[width=3.4cm,height=3.4cm]{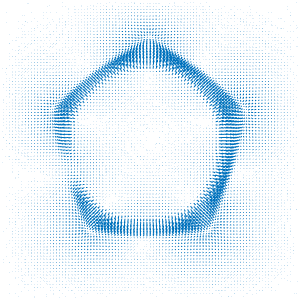}}
	\subfigure[$\pmb{D}_{1}$]{\includegraphics[width=3.4cm,height=3.4cm]{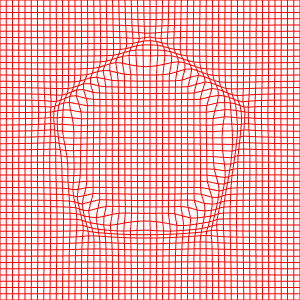}}
	\subfigure[$\pmb{D}_{1}\circ\pmb{C}$ vs $\pmb{P}$]{\includegraphics[width=3.4cm,height=3.4cm]{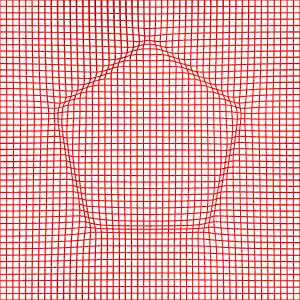}}	
	\subfigure[$\pmb{u}_{2}$]{\includegraphics[width=3.4cm,height=3.4cm]{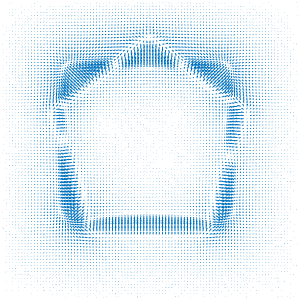}}
	\subfigure[$\pmb{D}_{2}$]{\includegraphics[width=3.4cm,height=3.4cm]{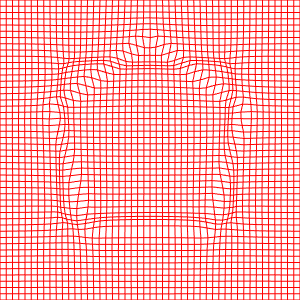}}
	\subfigure[$\pmb{D}_{2}\circ\pmb{P}$ vs $\pmb{R}$]{\includegraphics[width=3.4cm,height=3.4cm]{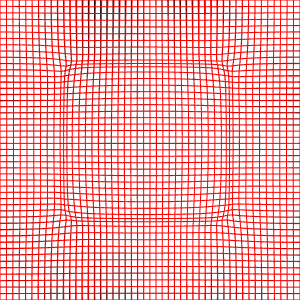}}
	\subfigure[$\pmb{D}=\pmb{D}_{2}\circ\pmb{D}_{1}$]{\includegraphics[width=3.4cm,height=3.4cm]{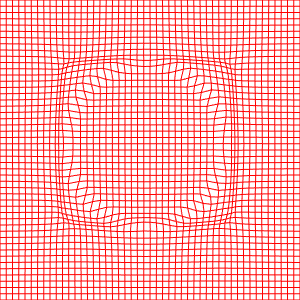}}
	\subfigure[$\pmb{D}\circ\pmb{C}$]{\includegraphics[width=3.4cm,height=3.4cm]{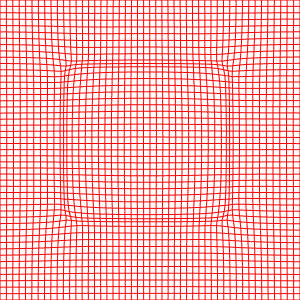}}

	\subfigure[$\pmb{u}_{3}$]{\includegraphics[width=3.4cm,height=3.4cm]{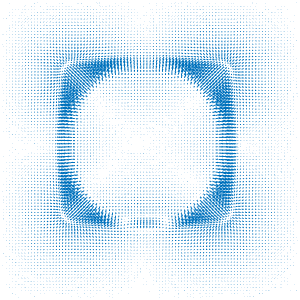}}
	\subfigure[$\pmb{D}_{3}$]{\includegraphics[width=3.4cm,height=3.4cm]{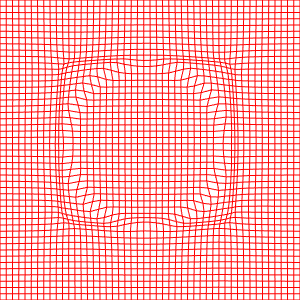}}
	\subfigure[$\pmb{D}_{3}\circ\pmb{C}$ vs $\pmb{R}$]{\includegraphics[width=3.4cm,height=3.4cm]{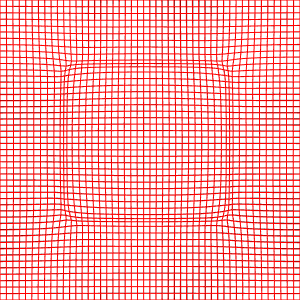}}
	\subfigure[$\pmb{D}_{3}$ vs $\pmb{D}$]{\includegraphics[width=3.4cm,height=3.4cm]{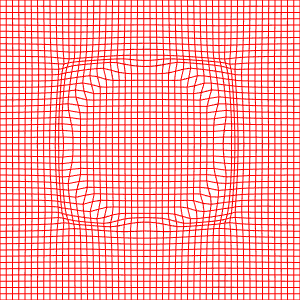}}
	\caption{Transitivity}\label{trans}
\end{figure}

\begin{table}[h]
	\small{
		\begin{center}
			\begin{tabular}{|c|c|c|c|ccc|}
				\hline	
				\multirow{2}{*}{orientation} & \multirow{2}{*}{$\mathrm{\Omega}$}& \multirow{2}{*}{$ratio$} & \multirow{2}{*}{sec $|$ iteration} &\multicolumn{3}{c|}{max differences of Fig.\ref{trans}(f, i, n)} \\  \cline{5-7} 
				&      &    &    &$|$det$\nabla(\_)|$&$|\nabla\times(\_)|$ &  $|| (\_) ||_{2}/|\mathrm{\Omega}|$     \\ \hline
				{\tt$\pmb{C}$ to $\pmb{P}$} &$[1,97]^2$ & $0.0013\%$  & 203.20 $|$ 16949  & 0.0349 &0.0222& $2.2357*10^{-4}$\\ 
				{\tt $\pmb{P}$ to $\pmb{R}$} &$[1,97]^2$ & $0.0020\%$ & 202.78 $|$ 19644   & 0.0399  &0.0326 &  $3.4142*10^{-5}$ \\
				{\tt$\pmb{C}$ to $\pmb{R}$} &$[1,97]^2$ & $0.0009\%$  & 193.76 $|$ 19495  & 0.0302 &0.0177& $1.4442*10^{-5}$ \\ 
				\hline
			\end{tabular}
			\caption{Performance of Fig.\ref{trans}}\label{tbl4}
		\end{center}
	}
\end{table}

\subsection{$\mathbf{Example}$: Find Inverse transformation of an Image Registration deformation}\label{eg35}
This example connects the revised VP to the proposed image registration method in \cite{Zhou} and provides a further insight about the solutions to the image registration method are potentially of the same diffeomorphism group that the revised VP focuses. An extensive and detailed study of this insight will be included in a different paper, while here is merely an intuitive demonstration. 

Image registration is the task of pixel value alignment that to move a moving image, $I_{J}$ in this example, to a fixed image, $I_{V}$ in this example, by a transformation, $\pmb{\phi}$, so that $I_{J}(\pmb{\phi})$ is similar to $I_{V}$. It is done with a minimization problem:
\begin{equation*}
	\min_{\phi} \int_{\mathrm{\Omega}} [I_{J}(\pmb{\phi})-I_{V}]^{2}.
\end{equation*} 
Ideally, a diffeomorphic $\pmb{\phi}$ is expected. In Fig.(\ref{J2V}), (a) is the moving image to be moved to the fixed image (b); (c) is the diffeomorphic solution found by \cite{Zhou}; (e) is the moved/registered image which is very similar to the fixed image (b).

\begin{figure}[H]
	\centering
	\subfigure[$I_{J}$]{\includegraphics[width=3.4cm,height=3.4cm]{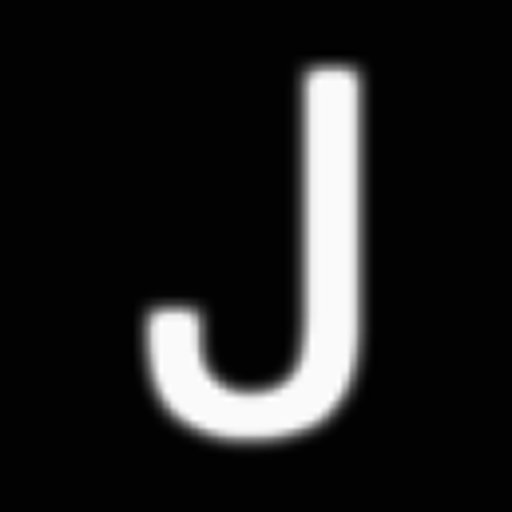}}
	\subfigure[$I_{V}$]{\includegraphics[width=3.4cm,height=3.4cm]{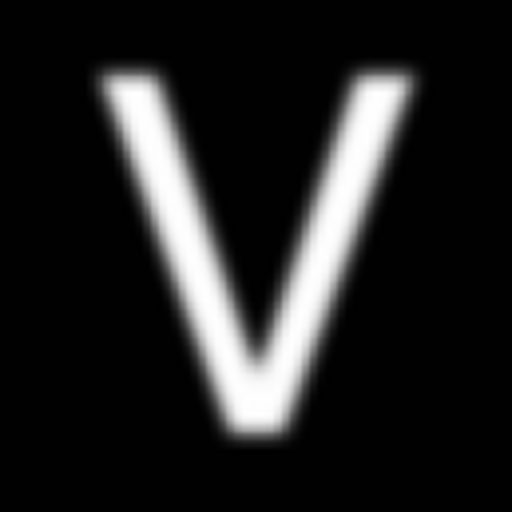}}

	\subfigure[$\pmb{\phi}$]{\includegraphics[width=3.4cm,height=3.4cm]{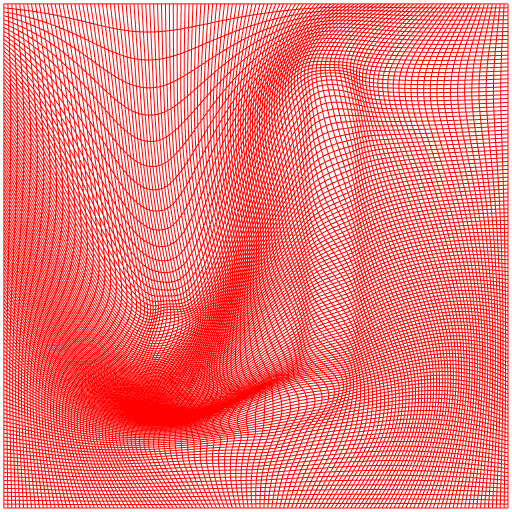}}
	\subfigure[$\pmb{u}_{\pmb{\phi}}$]{\includegraphics[width=3.4cm,height=3.4cm]{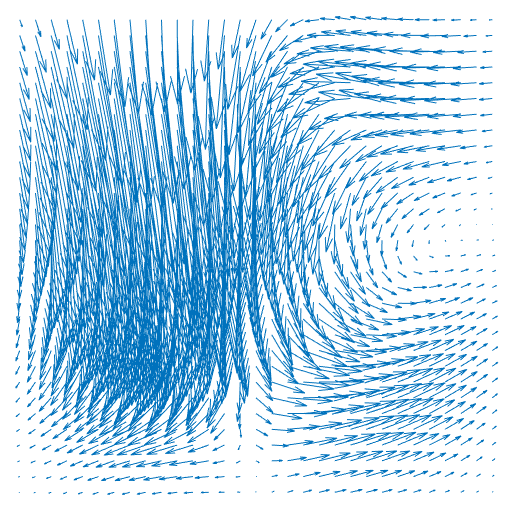}}
	\subfigure[$I_{J}(\pmb{\phi})$]{\includegraphics[width=3.4cm,height=3.4cm]{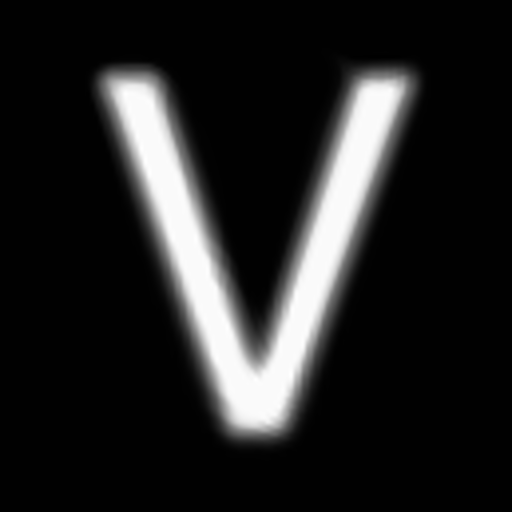}}
\end{figure}

Next in Fig.(\ref{J2V}), (f) is the guessed inverse transformation of (c) found by the revised VP; in (h), the composition of (f) and (c) in red grid is superposed on black grid $\pmb{id}$ but the black grid lines are barely shown, therefore,  (f) is a very good guess and is in the same diffeomorphism group where (c) is in. So, the interesting question is whether (f) is also the inverse transformation that moves the fixed image (b), $I_{V}$, back to the moving image (a), $I_{J}$, in the sense of image registration; (i) gives a positive response, i.e., (f) can be treated as the registration deformation from (b) to (a). 
\begin{figure}[H]
	\centering		
	\subfigure[$\pmb{\phi}^{-1}$ by revised VP]{\includegraphics[width=3.4cm,height=3.4cm]{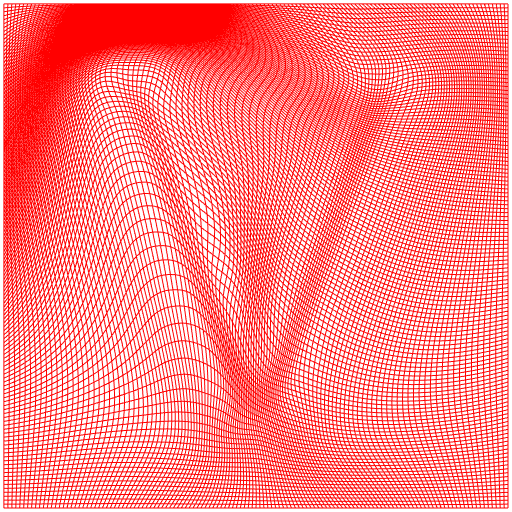}}
	\subfigure[$\pmb{u}_{\pmb{\phi}^{-1}}$]{\includegraphics[width=3.4cm,height=3.4cm]{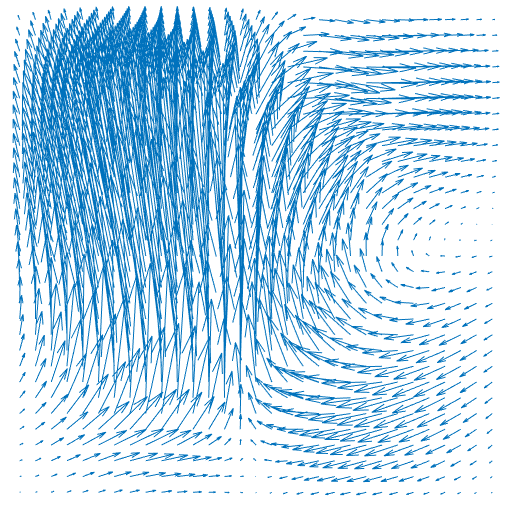}}
	\subfigure[$\pmb{\phi}^{-1}\circ\pmb{\phi}$ vs $\pmb{id}$]{\includegraphics[width=3.4cm,height=3.4cm]{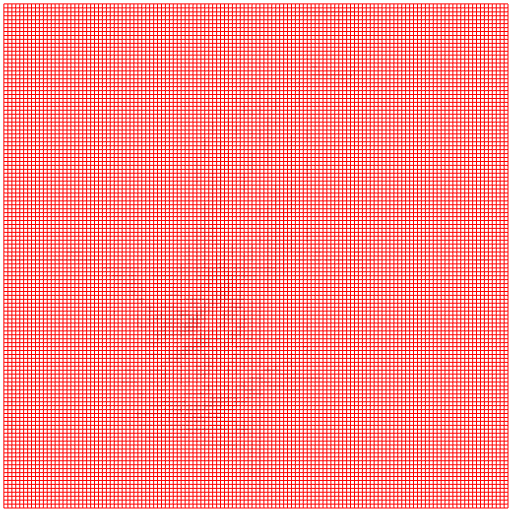}}
	\subfigure[$I_{V}(\pmb{\phi}^{-1})$]{\includegraphics[width=3.4cm,height=3.4cm]{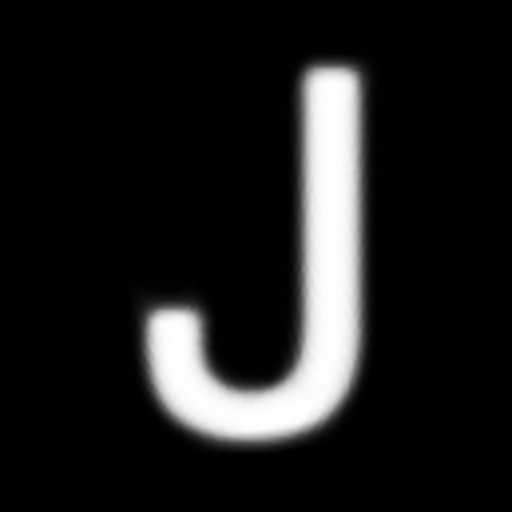}}
	\caption{Inverse deformation guessed by revised VP}\label{J2V}
\end{figure}

\section{Proposed Strategy to handle the Mismatch Issue with Preliminary Results}\label{demo2}
It is observed that optimizing along JD of (\ref{ssddelfg}) tends to converge fast in a global view and along curl of (\ref{ssddelfg}) tends to meet in local conditions slowly. And, alternating these two directions does not necessarily produce desirable results, we propose a computational strategy (proposed strategy) to bypass the mismatch issue from a matched point of view by taking advantages of revised VP. It goes as follows:
\begin{itemize}
	\item[$\bullet$] Step-1, reconstruct $\pmb{\phi}_{jd}$ with $f_{o}$, only along JD of (\ref{ssddelfg}), since it converge faster in more of a global manner, so $\pmb{\phi}_{jd}$ is expected in lacking information of $\pmb{g}_{o}$;

	\item[$\bullet$] Step-2, form $\pmb{\phi}=\pmb{\phi}_{c}(\pmb{\phi}_{jd})$ by reconstruct $\pmb{\phi}_{\pmb{m}}=\pmb{\phi}_{c}$ with $f_{o}$, $\pmb{g}_{o}$, and $\pmb{\phi}_{o}=\pmb{\phi}_{jd}$, which allows $SSD$ be optimized mainly to meet curl in a local manner along (\ref{ssd2delF}), even (\ref{ssd2delF}) represents the gradient in a combined way.
\end{itemize}
Both steps are achieved by the algorithm above with correct choice of inputs and proper codding of key inner lines. E.g.\ref{eg4} presents the work flow of the proposed strategy. E.g.\ref{eg5} shows that the proposed strategy and the revised VP catch noised curl better than the original VP. E.g.\ref{eg6} describes how they respond to leveled distortions from matched pairs of JD and curl to unmatched pairs.
\subsection{$\mathbf{Example}$: Work Flow of the Proposed Strategy}\label{eg4}
Fig.\ref{sep}(a) is given purposely with large curl realized by a cut-off counter-clockwise rotation, so the effectiveness of the proposed strategy can be visibly shown. Fig.\ref{sep}(b,c,d) corresponds to Step-1 and the first row of Table.\ref{tbl1}. Fig.\ref{sep}(b,c) shows that only JD had been mainly recovered, where $\pmb{\phi}_{jd}$ does not align well with $\pmb{G}$ in Fig.\ref{sep}(d). Fig.\ref{sep}(e,f,g,h) corresponds to Step-2 and the second row of Table.\ref{tbl1}. Fig.\ref{sep}(e,f,g,h) reflects the unmet curl is recovered after Step-2. Fig.\ref{sep}(f) is the intermediate transformation $\pmb{\phi}_{\pmb{m}}$ that left-translates $\pmb{\phi}_{jd}$ to $\pmb{\phi}$. Fig.\ref{sep}(h) shows the reconstructed grid lines up with quite close to Fig.\ref{sep}(a). Only tiny portion of black grids can be seen. Fig.\ref{sep}(b,e) respectively are the displacements based on JD and curl, which are acquired individually and made possible by the proposed strategy.
\begin{figure}[H]
	\centering
	\subfigure[$\pmb{G}$]{\includegraphics[width=3.4cm,height=3.4cm]{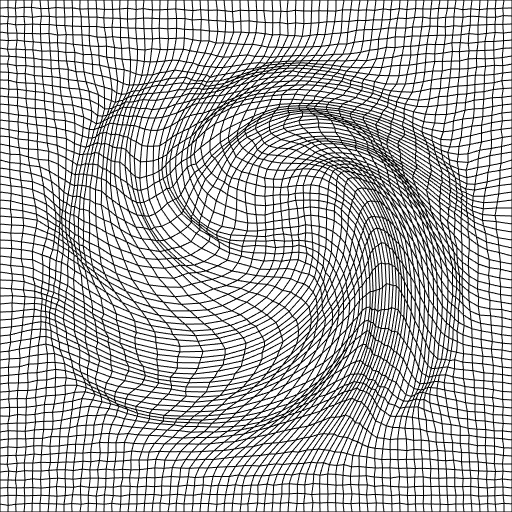}}
	\subfigure[$\pmb{u}_{jd}$]{\includegraphics[width=3.4cm,height=3.4cm]{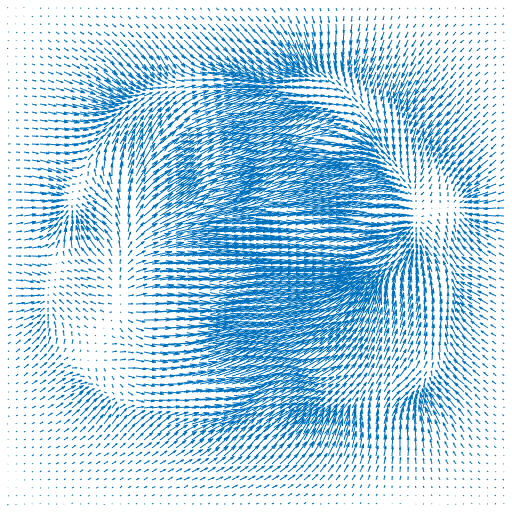}}
	\subfigure[$\pmb{\phi}_{jd}$]{\includegraphics[width=3.4cm,height=3.4cm]{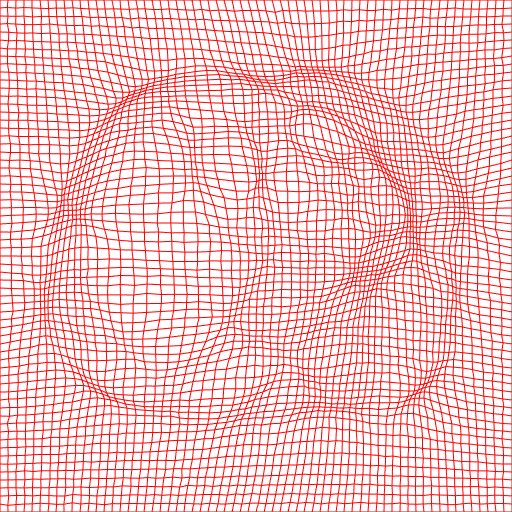}}
	\subfigure[$\pmb{\phi}_{jd}$ vs $\pmb{G}$]{\includegraphics[width=3.4cm,height=3.4cm]{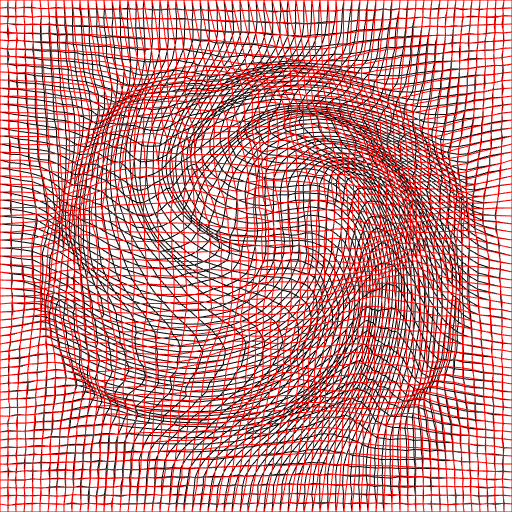}}

	\subfigure[$\pmb{u}_{c}$]{\includegraphics[width=3.4cm,height=3.4cm]{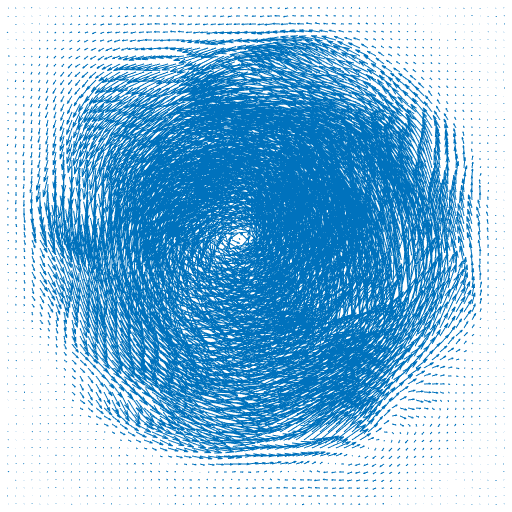}}
	\subfigure[$\pmb{\phi}_{c}$]{\includegraphics[width=3.4cm,height=3.4cm]{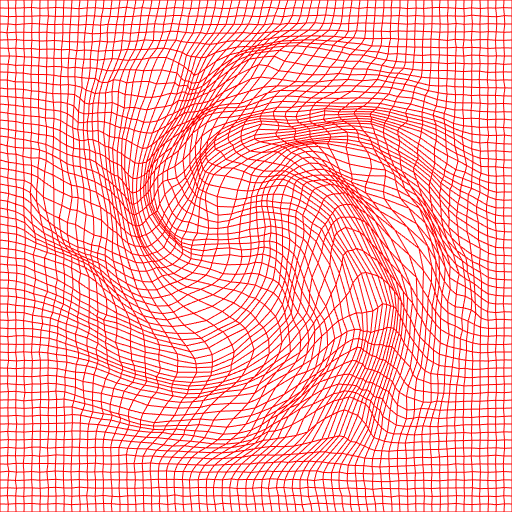}}
	\subfigure[$\pmb{\phi}=\pmb{\phi}_{c}\circ\pmb{\phi}_{jd}$]{\includegraphics[width=3.4cm,height=3.4cm]{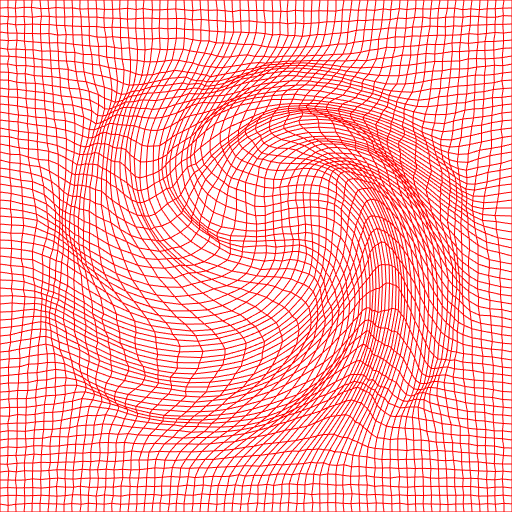}}
	\subfigure[$\pmb{\phi}$ vs  $\pmb{G}$]{\includegraphics[width=3.4cm,height=3.4cm]{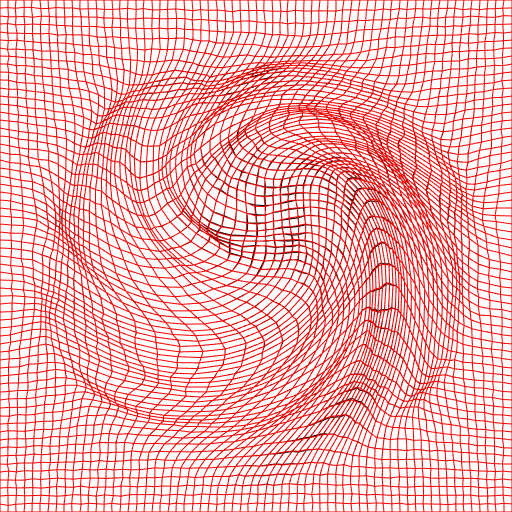}}
	\caption{Work Flow of the Proposed Strategy}\label{sep}
\end{figure}

\begin{table}[h]
	\small{
		\begin{center}
			\begin{tabular}{|c|c|c|c|ccc|}
				\hline	
				\multirow{2}{*}{orientation} & \multirow{2}{*}{$\mathrm{\Omega}$}& \multirow{2}{*}{$ratio$} & \multirow{2}{*}{sec$|$iteration} &\multicolumn{3}{c|}{max differences of  Fig.\ref{sep}(d, h) }  \\  \cline{5-7} 
				&      &    &    & $|$det$\nabla(\_)|$ & $|\nabla\times(\_)|$ &  $|| (\_) ||_{2}/|\mathrm{\Omega}|$     \\ \hline
				{\tt$\pmb{id}$ to $\pmb{\phi}_{jd}$} &$[1,65]^2$ & $81.138\%$& 15.04 $|$ 1779  & 0.1935 &2.6428& $1.4369*10^{-3}$ \\ 
				{\tt$\pmb{\phi}_{jd}$ to $\pmb{G}$} &$[1,65]^2$ & $0.0412\%$   & 118.61 $|$ 20000  & 0.1058 &0.0734& $2.7947*10^{-5}$ \\ 
				\hline
			\end{tabular}
			\caption{Performance of E.g.\ref{eg4}}\label{tbl1}
		\end{center}
	}
\end{table}
\subsection{$\mathbf{Example}$: Comparison in Catching the Noises}\label{eg5}
Since JD$>0$ is the key to maintain diffeomorphic solutions, so, to see how VP is affected from a noised curl, in this example, det$\nabla\pmb{G}=f_o$ and $\nabla\times\pmb{G}=\pmb{g}_{o}$ of E.g.\ref{eg4} are taken, but a Gaussian noise realized with $\pm67.5^{\circ}$ local rotation is added over a disk window to $\pmb{g}_{o}$, denoted as $\pmb{\hat{g}}_{o}$. Hence, there is no ground truth to compare with. Fig.\ref{sep2}(a,b,c,d,e) are solutions from the proposed strategy, Fig.\ref{sep2}(f) is of the revised VP and Fig.\ref{sep2}(g) is of the original VP. 
\begin{figure}[H]
	\centering
	\subfigure[$\pmb{\phi}_{jd}$]{\includegraphics[width=3.4cm,height=3.4cm]{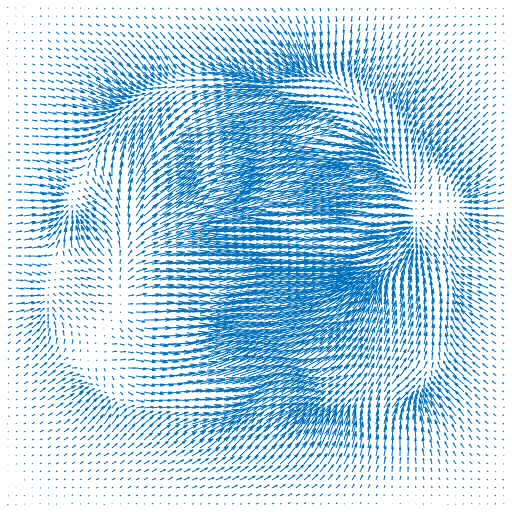}}
	\subfigure[$\pmb{u}_{jd}$]{\includegraphics[width=3.4cm,height=3.4cm]{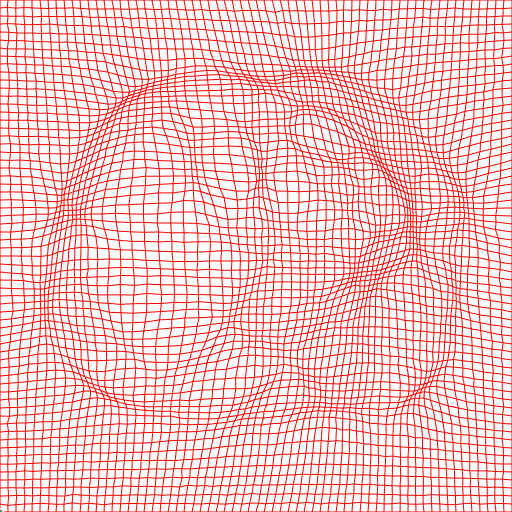}}
	\subfigure[$\pmb{u}_{c}$]{\includegraphics[width=3.4cm,height=3.4cm]{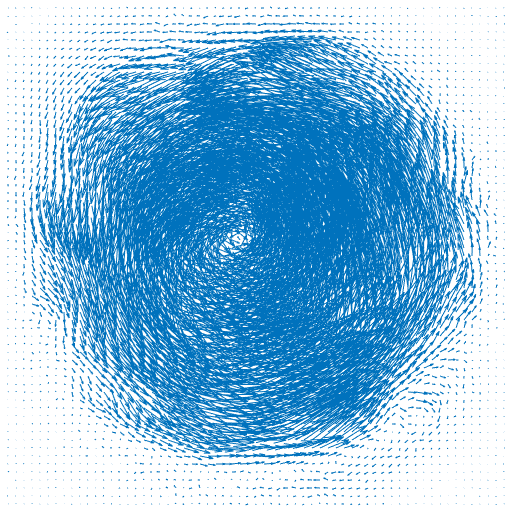}}

	\subfigure[$\pmb{\phi}_{c}$]{\includegraphics[width=3.4cm,height=3.4cm]{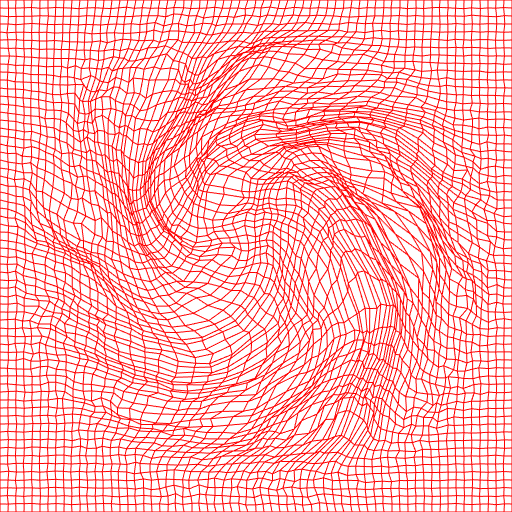}}
	\subfigure[$\pmb{\phi}=\pmb{\phi}_{c}\circ\pmb{\phi}_{jd}$]{\includegraphics[width=3.4cm,height=3.4cm]{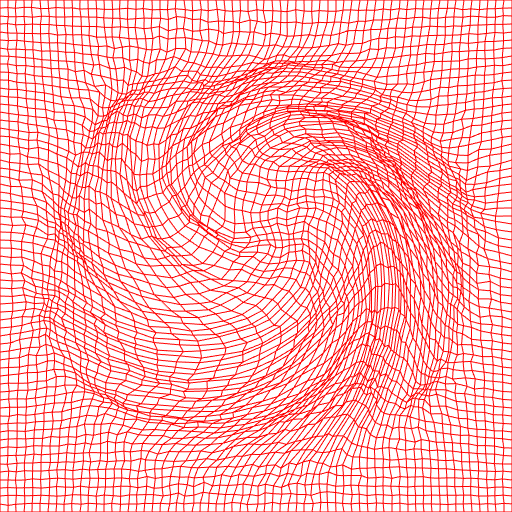}}
	\subfigure[$\pmb{\phi}_{rvp}$]{\includegraphics[width=3.4cm,height=3.4cm]{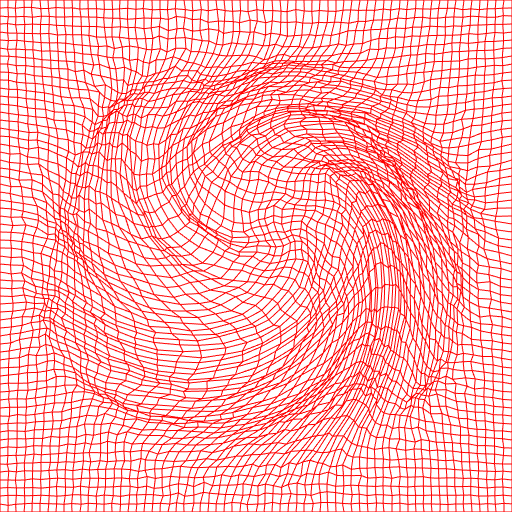}}
	\subfigure[$\pmb{\phi}_{ovp}$]{\includegraphics[width=3.4cm,height=3.4cm]{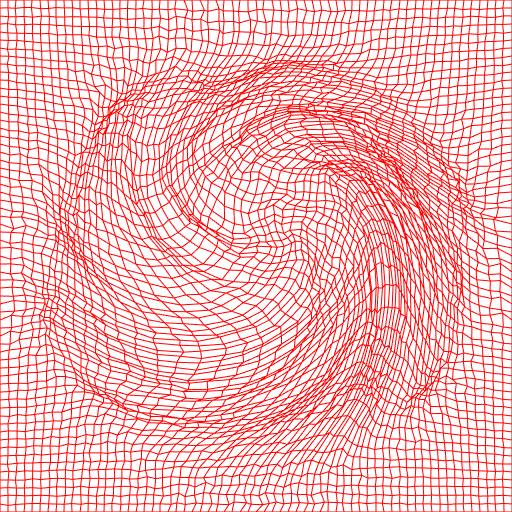}}
	\caption{New Version compare to Original VP}\label{sep2}
\end{figure}

\begin{table}[h]
	\small{
		\begin{center}
			\begin{tabular}{|c|c|c|c|c|c|c|}
				\hline	
				orientation & $ratio$ &max JD &min JD &max$|$det$\nabla(\_)-f_{o}|$ &max$|\nabla\times(\_)-\pmb{\hat{g}}_{o}|$    \\ \hline
				{\tt$\pmb{id}$ to $\pmb{\phi}_{jd}$}  to $\pmb{\phi}$& $0.5565\%$  & 2.3334  &0.1465    & 0.2041 &0.2649 \\ 
				{\tt$\pmb{id}$  to $\pmb{\phi}_{rvp}$} &$0.3103\%$  & 2.3357  & 0.1677  & 0.1976 &0.2367\\ 
				{\tt$\pmb{id}$  to $\pmb{\phi}_{ovp}$} &$0.2164\%$  & 2.3370  & 0.1714  & 0.3260 &0.7637\\ 
				\hline
			\end{tabular}
			\caption{Performance of E.g.\ref{eg5}}\label{tbl2}
		\end{center}
	}
\end{table}
In Table.\ref{tbl2}, the $ratio$'s indicate that all solutions converged when stopping criteria is met; the max and min JD values suggest that the proposed strategy, the revised VP and original VP all produce very comparable results in terms of smoothness; but, the max differences in JD and curl point that the proposed strategy and the revised VP are more truthful in satisfying local JD and curl than the original VP. Observed this, it should be noted that a lower value of $ratio$ does not necessarily mean a better solution in general, especially, when unmatched issue arises. 

\subsection{$\mathbf{Example}$: Comparison in Responding to Distortions}\label{eg6}
Take det$\nabla\pmb{G}=f_o$, $\nabla\times\pmb{G}=\pmb{g}_{o}$ of E.g.\ref{eg4}, but $\pmb{g}_{o}$ is distorted in from 0 to 10 levels of clockwise rotation over a 16-by-16 window. The window is located at north-eastern side where the grid-cells tend to entangle. This example stimulates the flow from a matched pair of JD and curl becoming an unmatched pair. Again, since JD$>0$ is the condition for producing diffeomorphisms, so the distortions are only applied to the prescribed curl. Similarly, there is no ground truth here. It is expected to see that the proposed strategy should handle higher level of distortions than the revised VP and the original VP. Fig.\ref{sep5}(a-h) are the results of proposed strategy in level 5-8, Fig.\ref{sep5}(i-p) are the result of revised VP in level 5-8 and Fig.\ref{sep5}(q-x) are the result of original VP in level 5-8. 
\begin{figure}[H]
	\centering
	\subfigure[$\pmb{\phi}_{ps}$ lvl-5]{\includegraphics[width=3.4cm,height=3.4cm]{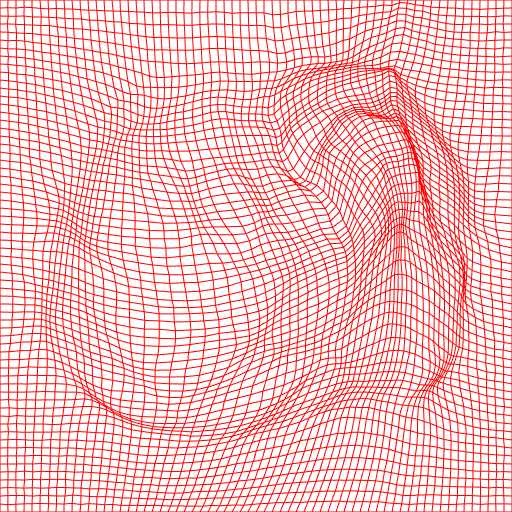}}
	\subfigure[$\pmb{\phi}_{ps}$ lvl-6]{\includegraphics[width=3.4cm,height=3.4cm]{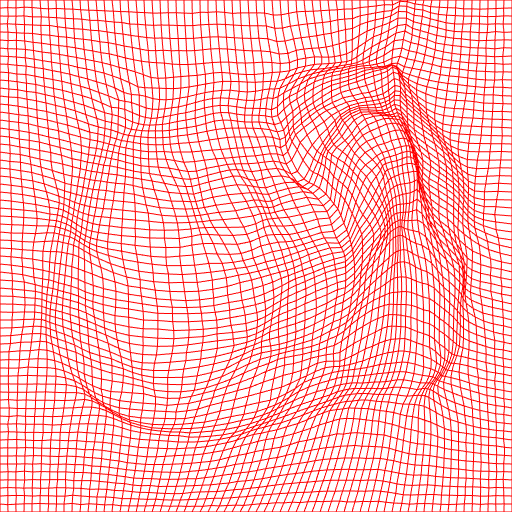}}
	\subfigure[$\pmb{\phi}_{ps}$ lvl-7]{\includegraphics[width=3.4cm,height=3.4cm]{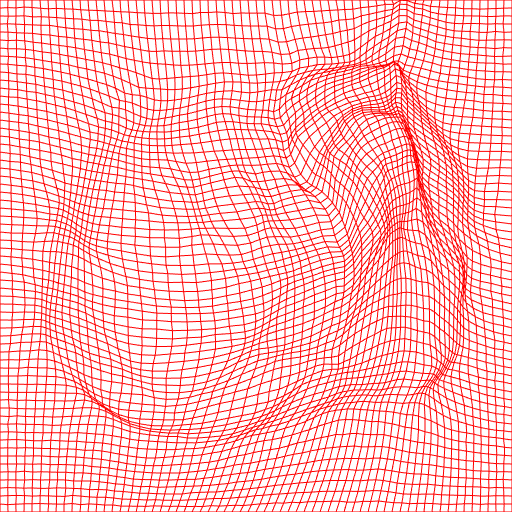}}
	\subfigure[$\pmb{\phi}_{ps}$ lvl-8]{\includegraphics[width=3.4cm,height=3.4cm]{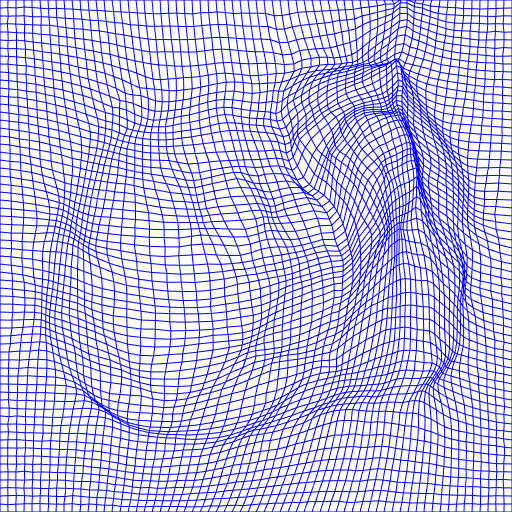}}
 
	\subfigure[$\pmb{\phi}_{rvp}$ lvl-5]{\includegraphics[width=3.4cm,height=3.4cm]{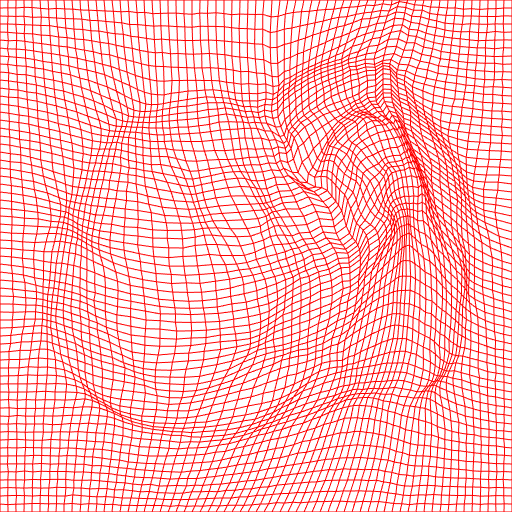}}
	\subfigure[$\pmb{\phi}_{rvp}$ lvl-6]{\includegraphics[width=3.4cm,height=3.4cm]{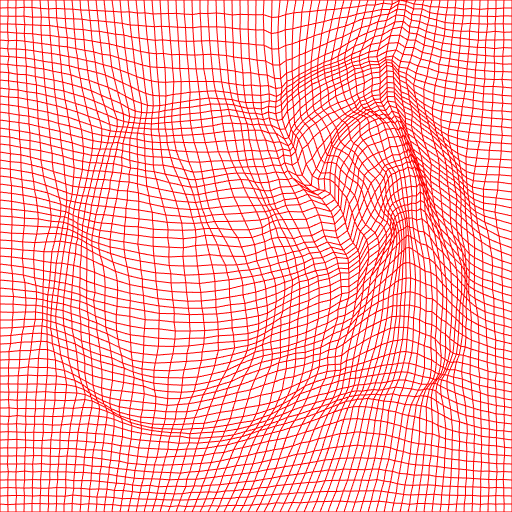}}
	\subfigure[$\pmb{\phi}_{rvp}$ lvl-7]{\includegraphics[width=3.4cm,height=3.4cm]{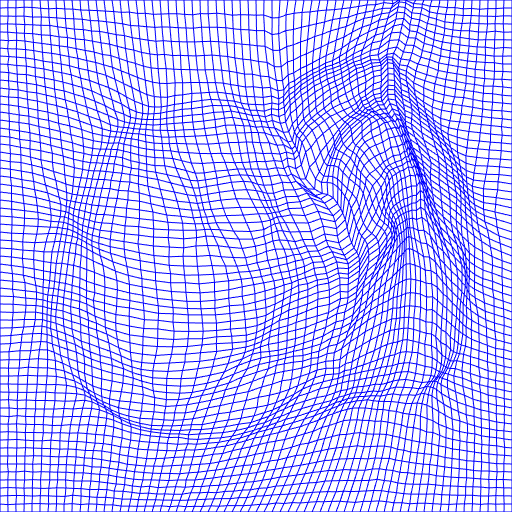}}
	\subfigure[$\pmb{\phi}_{rvp}$ lvl-8]{\includegraphics[width=3.4cm,height=3.4cm]{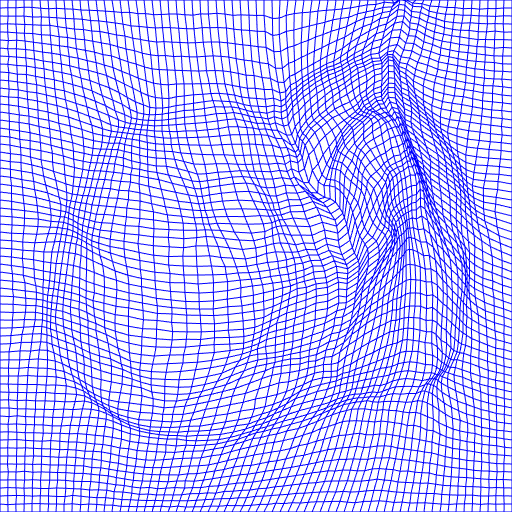}}
 
	\subfigure[$\pmb{\phi}_{ovp}$ lvl-5]{\includegraphics[width=3.4cm,height=3.4cm]{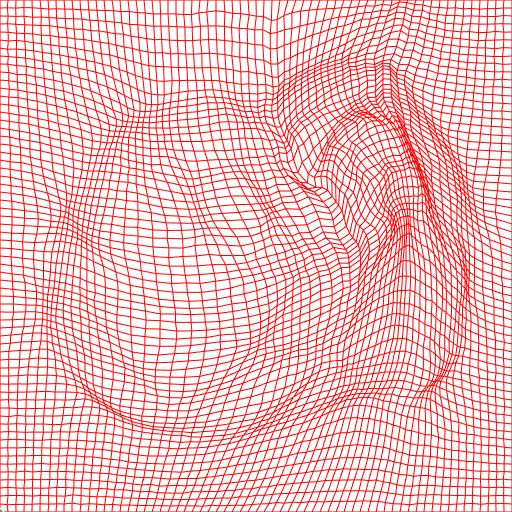}}
	\subfigure[$\pmb{\phi}_{ovp}$ lvl-6]{\includegraphics[width=3.4cm,height=3.4cm]{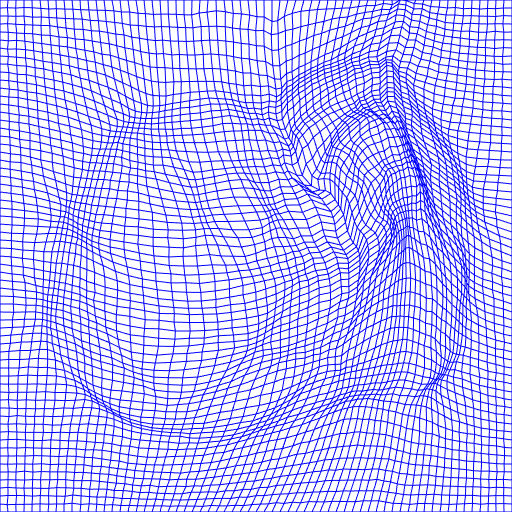}}
	\subfigure[$\pmb{\phi}_{ovp}$ lvl-7]{\includegraphics[width=3.4cm,height=3.4cm]{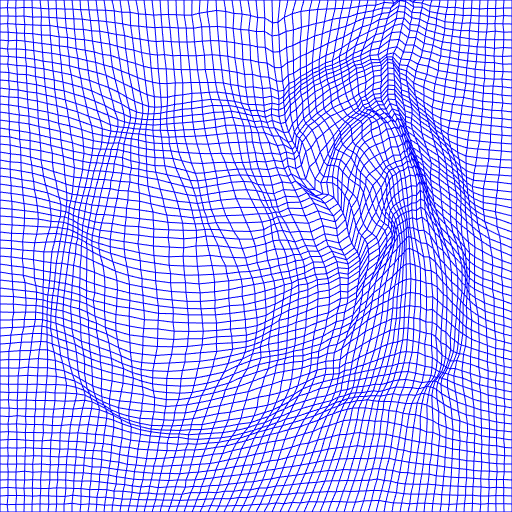}}
	\subfigure[$\pmb{\phi}_{ovp}$ lvl-8]{\includegraphics[width=3.4cm,height=3.4cm]{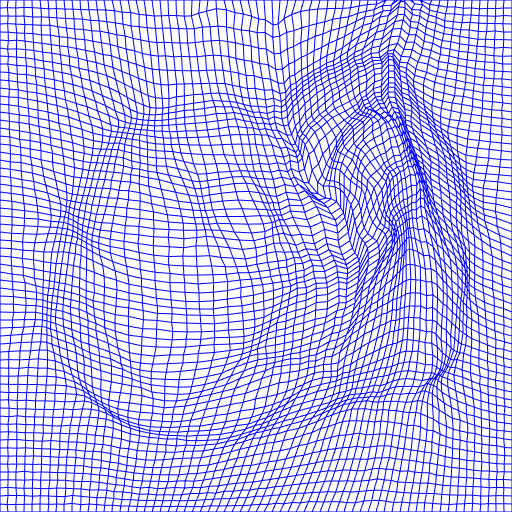}}
	\caption{Comparison in Distortion level 5 to 8}\label{sep5}
\end{figure}
Fig.\ref{sep6}(a) records all min JDs of the solutions, where the boxed values point out at which level that min JD turned negative. It reflects that as the distortion level increases, the proposed strategy maintain diffeomorphic solutions up to level 7 while the revised VP only stay diffeomorphic to level 6 and the original VP failed to surpass level 5. In Fig.\ref{sep5}, the diffeomorphic solutions are plotted in red while the non-diffeomorphic solutions are plotted with blue grids, due to the challenge of visualizing their differences in tiny cells.
\begin{figure}[H]
	\centering
	\subfigure[min JD v.s. lvl]{\includegraphics[width=4.0cm,height=3.0cm]{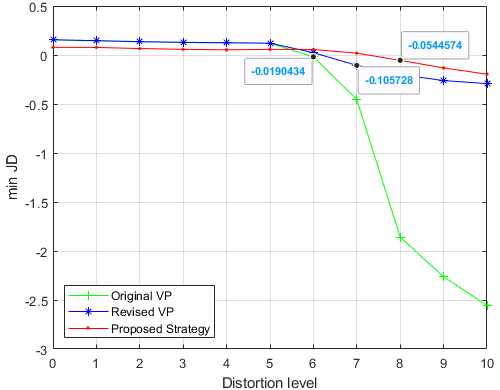}}
	\subfigure[max JD v.s. lvl]{\includegraphics[width=4.0cm,height=3.0cm]{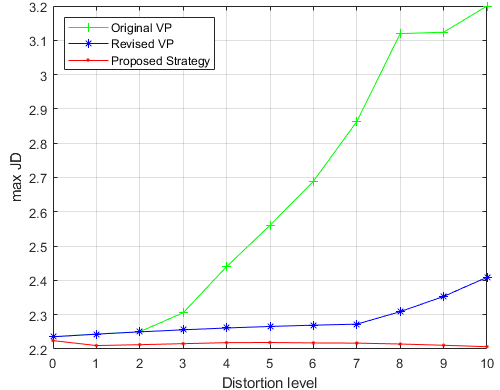}}
	\subfigure[$ratio$ v.s. lvl]{\includegraphics[width=4.0cm,height=3.0cm]{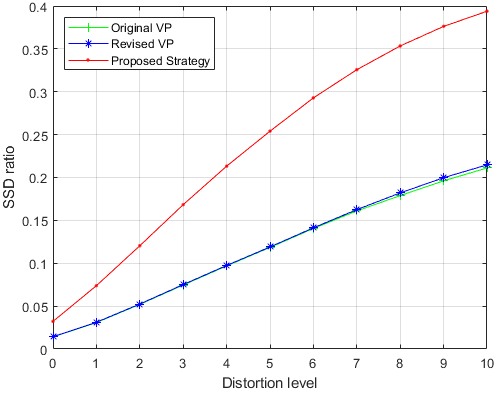}}
\end{figure}

\begin{figure}[H]
\centering
	\subfigure[max-diff JD v.s. lvl]{\includegraphics[width=4.0cm,height=3.0cm]{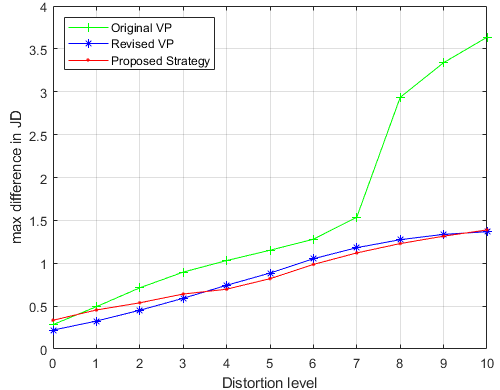}}
	\subfigure[max-diff curl v.s. lvl]{\includegraphics[width=4.0cm,height=3.0cm]{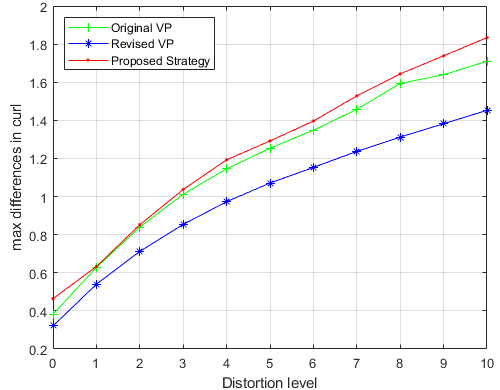}}
	\caption{Performance Graphs of Fig.\ref{sep5}}\label{sep6}
\end{figure}
Fig.\ref{sep6}(a, b) recorded the proposed strategy and the revised VP tend to maintain smaller gaps between max and min JD. This means the proposed strategy and revised VP produce smoother solutions compare to original VP. Fig.\ref{sep6}(c) shows that the revised VP and original VP shrink $SSD$ regardless JD and curl are matched or not. Thanks to Step-1 of the proposed strategy, its solutions are more loyal in meeting $f_o$ rather simply reducing $SSD$. In Fig.\ref{sep6}(d, e), both the proposed strategy and revised VP presented more potentials to stay close to $f_o$ when curl becomes unmatched to $f_o$. This example suggest that proposed strategy has better chance to handle an unmatched scenario and the revised VP is more stable in a matched situation. Overall, the proposed strategy and the revised VP have their own advantages over the original VP. 

\section{Discussion}
The original VP is unsuccessful in characterizing of its solutions are of the diffeomorphism group that collects from $H^{2}_{0}(\mathrm{\Omega})$ and has poor performance in 3D grid generations. These limitations left the original VP acceptable in the realm of grid generation for certain engineering purposes, but once it comes to the areas that demand higher accuracy such as medical image processing, undesired results can be accumulated. Actually, this revision of VP also aims at providing a better theoretical insight and computation tool to enhance the novel image atlas construction method proposed in \cite{Zhou}, as it shown in E.g.(\ref{eg35}). The studies in this regard will be included in future papers.

\section{Conclusion}
In sum, the revised VP is enabled to construct transformations that cope with composition of transformations as the left-translation of the diffeomorphism group and welcomes more candidates to be characterized by using only JD and curl. It handles 3D grid generations much better than the original VP. Another contribution of this paper is the 2-steps new computational strategy which allows JD and curl information to be acquired in a separately manner. It is also built as a preliminary tool to handle the mismatched issue. Analytical study of how the provided algorithm converge is of our interests. For future, we will build more efficient algorithms as preparation for 3D real medical image problems.

\vspace{-0.3cm}

\singlespacing


\begin{thebibliography}{20}
\bibitem{Bauer}M.~Bauer, S. Joshi and K.~Modin, Diffeomorphic density matching by optimal information transport, \emph{SIAM Journal on Imaging Sciences},\textbf{8}:3 (2015), 1718–1751.	
\bibitem{Brackbill}J.U.~Brackbill and J.S.~Saltzman, Adaptive zoning for singular problems in two dimensions, \emph{J. Comput. Phys},\textbf{46}:(1982), 342-368.	
\bibitem{Cai} X.~Cai, D.~Fleitas, B.~Jiang and G.~Liao, Adaptive grid generation based on the least-squares finite elements method,\emph{Computers and Mathematics with Applications},\textbf{48}:(2004),1007-1085.	
\bibitem{Castillo}J.~Castillo, S.~Steinberg, and P.~Roach, Parameter estimation in variational grid generation,\emph{Appl. Math. Comput.},\textbf{155}:(1988), 155-177.	
\bibitem {ChenY}Y.~Chen and X.~Ye, Inverse Consistent Deformable Image Registraition.
\emph{The Legacy of Alladi Ramakrishnan in the Mathematical Sciences, SpringerVerlag}, (2010,) 419-440.
\bibitem{ChenXi}X.~Chen and G.~Liao, New Variational Method of Grid Generation with Prescribed Jacobian determinant and Prescribed Curl. preprint, arxiv.org/pdf/1507.03715, 2015.
\bibitem{DacMos} B.~Dacoragna and J.~Moser, On A Partial Differential Equation Involving the Jacobian determinant, \emph{Ann.Inst H Poincaré},\textbf{7}:(1990), 1-26.	
\bibitem{Gu}X.~Gu and S.T.~Yau, \emph{Computational Conformal Geometry}, International Press of Boston, Inc, 2008.
\bibitem{Grajewski}
M.~Grajewski, M.~Koster and S.~Turek, Mathematical and Numerical Analysis of a Robust and Efficient Grid Deformation Method in the Finite Element Context, \emph{SIAM Journal on Scientific Computing},  \textbf{31}:2 (2009), 1539-1557.
\bibitem{Huang} W.~Huang and W.~Sun, Variational mesh adaptation II: Error estimates and monitor functions, \emph{J. Comput. Phys}, \textbf{184}:(2003), 619-648.
\bibitem{Lee}E.~Lee and M.~Gunzburger, An Optimal Control Formulation of an Image Registration Problem, \emph{Journal of Mathematical Imaging and Vision}, \textbf{36}:1 (2009), 69-80.
\bibitem{Joshi} S,~Joshi, B.~Davis, M,~Jomier and G.~Gerig, Unbiased Diffeomorphic Atlas Construction for Computational Anatomy,\emph{NeuroImage}, \textbf{23}:(2004), 151-160.
\bibitem{Joshi2} S,~Joshi and M.~Miller, Landmark Matching via Large Deformation Diffeomorphisms, \emph{IEEE Transactions on Medical Imaging}, \textbf{9}:8 (2000), 1357-1370.
\bibitem{Liao} G.~Liao, X.~Cai, J.~Liu, X.~Luo, J.~Wang, J.~Xue, Construction of differentiable transformations, \emph{Applied Math.Letters},\textbf{22}:(2009), 543-1548.	
\bibitem {Liseikin}V.~Liseikin, \emph{Grid Generation Method},Springer Press. 1999.
\bibitem{Moser} J.~Moser, \emph{Volume elements of a Riemann manifold, Trans. AMS}, \textbf{120}: (1965), 155-177.	
\bibitem {Sotiras} A.~Sotiras, C.~Davatzikos and N,~Paragios, Deformable Medical Image Registration: A Survey, \emph{IEEE Transactions on Medical Imaging}, \textbf{32}: 7 (2013), 155-177.
\bibitem{Zhou} Z.~Zhou, \emph{Image Analysis Based on Differential Operators with Applications to Brain MRIs}, PhD Dissertation, University of Texas at Arlington, 2019.

\end{thebibliography}

\end{document}